\documentclass[11pt,b5paper,twoside,headrule]{amsart}

\usepackage{latexsym}
\usepackage{amssymb}

 \usepackage{amscd}

\footskip=10mm\parskip 0.2cm \addtocounter{page}{836}
\begin{document}

\numberwithin{equation}{section}
 \newtheorem{theorem}{Theorem}[section]
 \newtheorem{corollary}[theorem]{Corollary}
 \newtheorem{conj}[theorem]{Conjecture}
 \newtheorem{prop}[theorem]{Proposition}
 \newtheorem{lemma}[theorem]{Lemma}
 \newtheorem{remark}[theorem]{Remark}
 \newtheorem{axiom}[theorem]{Axiom}
 \newtheorem{defn}[theorem]{Definition}
 \newtheorem{examp}[theorem]{Example}
 \newtheorem{proposition}[theorem]{Proposition}
 \renewcommand{\proof}{{\noindent{\bf Proof.}\hspace{5pt}}}

\newcommand{\quott}{/\! /}

\newcommand{\liet}{{\bf t}}
\newcommand{\liets}{ { \bf t}^*}
\newcommand{\lieg}{{\bf g}}

\newcommand{\tri}{\bigtriangleup}
\renewcommand{\Delta}{\tri}

 \newcommand {\quoff}{/\!/}
\font\tenhtxt=eufm10 scaled \magstep0 \font\tenBbb=msbm10 scaled
\magstep0 \font\tenrm=cmr10 scaled \magstep0 \font\tenbf=cmb10
scaled \magstep0


\def\evenhead{{\protect\centerline{\textsl{\large{Lisa Jeffrey and Joon-Hyeok Song}}}\hfill}}

\def\oddhead{{\protect\centerline{\textsl{\large{Intersection Numbers in Reduced Spaces}}}\hfill}}

\pagestyle{myheadings} \markboth{\evenhead}{\oddhead}
\thispagestyle{empty}
\noindent{{\small\rm Pure and Applied Mathematics Quarterly\\ Volume 2, Number 3\\
(\textit{Special Issue: In honor  of \\ Robert MacPherson, Part 1 of 3})\\
837---865, 2006} \vspace*{1.5cm} \normalsize

\begin{center}
{\bf{\Large Intersection Numbers in Quasi-Hamiltonian\\ Reduced
Spaces}}
\end{center}

\begin{center}
\large{Lisa Jeffrey and Joon-Hyeok Song}

\it Dedicated to Bob MacPherson on the occasion of his 60th birthday
\end{center}

\footnotetext{Received July 17, 2005.}

\begin{center}
\begin{minipage}{5in}
\noindent \textbf{Abstract:} Jeffrey and Kirwan \cite{JK3} gave
expressions for intersection pairings on the reduced space of a
particular Hamiltonian $G$-space in terms of iterated residues. The
 definition of quasi-Hamiltonian spaces was introduced in \cite{AMM}. In
\cite{AMW1} a localization formula
 for  equivariant de Rham
cohomology of a compact quasi-Hamiltonian $G$-space was proved.
 In this paper we
 prove a residue formula for intersection pairings of  reduced spaces of
 quasi-Hamiltonian $G$-spaces,
by constructing a corresponding Hamiltonian $G$-space.
Our formula is a close analogue of   the result in \cite{AMW1}.
In this article  we rely heavily on the methods of \cite{JK3};
for the  general class of
  compact Lie groups $G$ treated in \cite{AMW1}, we rely
on results of Szenes and Brion-Vergne concerning diagonal bases.
\end{minipage}
\end{center}

\section{Introduction}
Alekseev, Meinrenken and
Woodward  \cite{AMW1}  (see also
\cite{AMW2}) proved a formula
for intersection numbers in quasi-Hamiltonian $G$-spaces, spaces which
had been  introduced
in \cite{AMM}. Our objective is to prove a close analogue
of the formula in \cite{AMW1}
 modelled on the proof for intersection
numbers in moduli spaces of representations of
the fundamental group of a Riemann surface \cite{JK3}, using the
residue theorem.

Our formula
 reduces to a sum over components
 $F$  of the fixed point set
 of subtori $S$ of $T$ for which $T/S$ acts locally
freely on $F$.
Our formula is a close analogue of  the formula in \cite{AMW1};
in Section 6 we comment on the differences between the two formulas.

We remark that the moduli spaces of interest in
\cite{JK3}
 are reduced spaces of a quasi-Hamiltonian space,
namely the product $G^{2h} $ where
$G= SU(n)$  and $h$ is the genus of the
Riemann surface.

The most important new results treated in our article are
as follows. First, we observe that the method used in
\cite{JK3} (by the
first author and F. Kirwan)
applies to quasi-Hamiltonian spaces more generally, not only to
$G^{2h} $. Second, we adapt this method using the
diagonal bases of Szenes and Brion-Vergne \cite{BV1,BV2,Sz} which
enables us to treat quasi-Hamiltonian spaces with an action of any
compact Lie group $G$, not only $G = SU(n)$.
Finally, our inductive treatment of the fixed point set (Theorem
\ref{t:residgen})
means that
it is unnecessary to identify the components of the fixed point set
at each stage of the induction.

The layout of this article is as follows. In Section 1 we recall
the properties of quasi-Hamiltonian $G$-spaces (spaces with group-valued
moment maps) and we state Szenes' theorem (which enables one to pass
between sums over the weight lattice and iterated residues).
Finally we recall the residue theorem which expresses intersection
numbers of reduced
spaces of Hamiltonian systems in terms of fixed point data on the original
Hamiltonian system.
In Section 2 we construct  the Hamiltonian  space corresponding to
a given quasi-Hamiltonian space.
In Section 3 we summarize material from \cite{JK3} on
equivariant Poincar\'e duals, and in section 4 we
summarize material from \cite{JK3} on
periodicity (an important property of Hamiltonian spaces corresponding
to quasi-Hamiltonian spaces).

 In Section 5 we
 prove a residue formula (Theorem \ref{t:residgen})
for intersection numbers
in the reduced space of a quasi-Hamiltonian $G$-space
(where $G$ is a compact Lie group).
In Section 6 we recall the results of Alekseev, Meinrenken and
Woodward \cite{AMW1} and give a proof of a close analogue of
the main  result of
\cite{AMW1}
by combining Szenes' theorem with our residue formula from Section 5.
We describe the differences between the two formulas.
Finally in Section 7 we describe some concrete examples.

The main result of \cite{AMW1} (Theorem 5.2 of that article)
is stated in terms of a sum over the dominant
weights $\lambda $ of a group $G$ (which are
identified with elements of $\liet$ using an
inner product on the Lie algebra of $G$)  and over the components $F$ of
the fixed point set of the vector field $v_\lambda$ on $M$ associated to
the action of $\lambda$.
Theorem 5.2 of \cite{AMW1} is closely related to Theorem \ref{t:residgen},
as we show using Szenes' theorem.
We need to use a very special case of Theorem 5.2 of \cite{AMW1}
(Proposition \ref{p:amwspcase}), the case when $G$ is a torus
acting locally freely on $M$; we expect that this special case may have a
proof independent of the proof given in \cite{AMW1}.

Our proof of  an analogue of Theorem 5.2 of \cite{AMW1}
 sheds light on
a novel  feature of quasi-Hamiltonian spaces, namely the
components  $F$ of the fixed point set of  some subtorus $S$ of $T$
for which $T/S$ acts locally freely on $F$.
Unlike the case of Hamiltonian $T$-actions, it is possible
that the fixed point set $M^T$ is empty although some connected subgroups
$S$ of $T$ have nonempty fixed point sets -- an example due to
Chris Woodward (treated in Section 7 (Example 7.3))
demonstrates how this situation may
arise.

  \subsection{Group-valued moment maps}

 Throughout this paper $G$ denotes a compact Lie group  of rank $r$
with Lie algebra $\bf g$.
A $G$-manifold is a compact manifold $M$
 together with a right action of $G$.
\begin{defn} \label{e:nudef}
For $\xi \in \bf{g}$ we denote by $\nu_\xi$
the  vector field on $M$ generated by $\xi$.
\end{defn}
 Let $\chi \in \Omega ^3 (G)$
denote the canonical closed bi-invariant 3-form on G:
\[
    \chi = \frac {1}{12} \langle \theta, [\theta, \theta]\rangle
 = \frac {1}{12} \langle\bar{\theta}, [\bar{\theta}, \bar{\theta}]\rangle.
\]
Here we have introduced a bi-invariant inner product
$\langle \cdot,\cdot\rangle$ on $\bf{g}$, and we denote by $\theta$ the
left-invariant Maurer-Cartan form, while $\bar{\theta}$
denotes the right-invariant Maurer-Cartan form.

 \begin{defn} {\cite{AMM}}
    A quasi-Hamiltonian $G$-space
is a $G$-manifold $M$ together with an invariant 2-form
   $\omega \in \Omega (M)^G$ and an equivariant map $\Phi \in C^{\infty}(M, G)^G$ (the quasi-Hamiltonian moment map) such that:
  \end{defn}

 (1) The differential of $\omega$ is given by:

    \[d\omega = - {\Phi}^* \chi\]

 (2) The map $\Phi$ satisfies

    \[\iota (\nu _\xi)\omega = \frac{1}{2} \Phi ^*
\langle \theta + \bar{\theta}, \xi\rangle \]

 (3) At each $x \in M$, the kernel of $\omega _x$  is given by

    \[\ker \omega _x = \{ \nu _{\xi}(x), \xi \in \ker(Ad_{\Phi(x)} + 1) \}\]

 \medskip

By analogy with
 Meyer-Marsden-Weinstein reduced spaces,
  quasi-Hamiltonian reduced spaces are defined
in the following way. Let $g \in G$ be a regular value of the
quasi-Hamiltonian  moment map $\Phi$. The preimage
$\Phi^{-1}(g)$ is a smooth manifold on which the action of the centralizer $G_g$ is locally free.
Then the reduced space $M_g = \Phi^{-1}(g)/G_g$ is a symplectic orbifold.

\begin{theorem} \label{t:qhr} (quasi-Hamiltonian reduction, \cite{AMM})
Let $M$ be a quasi-Hamiltonian $G_1 \times G_2$-space and $g \in G_1$ be a
regular value of the moment map $\Phi_1
:M \rightarrow G_1$. Then the pull-back of the 2-form $\omega$ to $\Phi_1^{-1}(g)$ descends to the reduced space
\[
    M_g = \Phi_1^{-1}(g)/G_g
\]
and makes it into a quasi-Hamiltonian $G_2$-space. In particular, if
$G_2 = \{e\}$ is trivial, then $M_g$ is a symplectic orbifold.
\end{theorem}

 \subsection{Szenes' Theorem}

\newcommand{\CC}{{\mathbb C}}
\newcommand{\RR}{{\mathbb R}}

\newcommand{\hypn}{{\mathcal H}}

\subsubsection{General $G$}
Let $G$ be a Lie group of rank $r$ with maximal torus $T$.
The general version of Szenes' theorem (Theorem \ref{t:szgen} below)
is proved in
\cite{Sz}; a different proof is given in \cite{BV2} (Theorem 27).
We shall require the following definitions.

The fundamental Weyl chamber is
denoted  ${\bf t}_+$; it  is a fundamental domain for the action
  of the Weyl group on ${\bf{t}}$.

We define variables
$\{Y_j , ~j = 1, \dots, r\}$
by $Y_j = e_j(X)$ for $X \in {\bf t}$ and
$e_j$ denoting the simple roots of $G$.

  If $g(Y_k,...,Y_{r})$ is a meromorphic function of
$Y_k,...,Y_{r}$, we interpret \linebreak$
Res_{Y_k=0}g(Y_k,...,Y_{r}) $ as the ordinary one-variable residue
of $g$ regarded as a function
  of $Y_k$ with $Y_{k+1},...,Y_{r}$ held constant.

\begin{defn} (\cite{BV2}) Let ${\bf V}$ be a real vector space and
let $\Delta \in {\bf V}^*$ be a finite collection of dual vectors.
The set
$$ \hypn  = \cup_{\alpha \in \Delta} \{ \alpha = 0 \} $$
is called a hyperplane arrangement.
An element in ${\bf V}$ is called {\em regular} if it is not in $\hypn$.
\end{defn}
\begin{remark} For our purposes ${\bf V} = \liet$, the Lie algebra of the
maximal torus $T$.
\end{remark}

\begin{defn}
 The ring $R_\Delta$ is the set of
rational functions with poles on
$\hypn$.
\end{defn}

\newcommand{\bases}{{\mathcal B}(\Delta)}

A subset $\sigma$ of $\Delta$ is called a {\em basis} of $\Delta$ if
the elements $\alpha \in \sigma$ form a basis of ${\bf V}$.
We denote the set of bases of $\Delta$ by $\bases$.
An {\em ordered basis} is a sequence of elements of $\Delta$ whose
underlying set is a basis.
\newcommand{\simplefrac}{S_\tri}

\begin{defn}
If $\sigma$ is a basis of $\Delta$ we let
$$\phi_\sigma (X) = \frac{1}{\prod_{\alpha \in \sigma} \alpha(X) } $$
and call $\phi_\sigma$ a {\em simple fraction}.
The collection of simple fractions is denoted
$\simplefrac$.
\end{defn}

\begin{defn} (\cite{BV2}, before
Definition 4) When $\tau = (\alpha_{i_1}, \dots, \alpha_{i_r} ) $ is
an ordered  basis of $\liet$ and
$X \in \liet \otimes \CC$, , we define variables
$Z_j^\tau = \alpha_{i_j} (X)$ for $j = 1, \dots, r$.
If $f$ is a meromorphic function on $\liet \otimes \CC$,
then we define
$$ {\rm Res}^{\tau} (f) =
{\rm Res}_{Z_1^\tau= 0 }  \dots   {\rm Res}_{Z_r^\tau= 0 }
(f)$$
for $X \in \liet \otimes \CC$.
\end{defn}

\begin{remark}
Note that the value of ${\rm Res}^\tau (f)$  depends on the
ordering of the basis $\tau$.
\end{remark}

\begin{defn} \label{d:diagbas} Let ${\bf V}$ be a real vector space.
A {\em diagonal basis} $OB$ is a subset of the ordered bases
${\mathcal OB} (\bigtriangleup)$ formed from  a collection
$$\bigtriangleup =\{\alpha_1, \dots,
\alpha_N \} $$ of elements   of ${\bf V}^*$
 (so that the null spaces of the
$\alpha_j$ specify
a collection of hyperplanes in ${\bf V}$)  satisfying the following
conditions:
\begin{enumerate}
\item  The set of simple fractions
$\phi_\sigma$ corresponding
to the  bases $\sigma \in OB$ forms a basis of $\simplefrac$
\item The collection
of meromorphic functions  on ${\bf V} \otimes \CC$ given by
$$ \phi_\sigma (X) =
\frac{1}{\alpha_{i_1}^\sigma(X) \dots\alpha^\sigma_{i_r}(X) }  $$
(for $\sigma \in { OB} $ given
by $\sigma =
(\alpha_{i_1}^\sigma, \dots,\alpha^\sigma_{i_r}) $)
satisfies
$${\rm Res}^{\tau}  (\phi_\sigma) =
\delta_\sigma^\tau,$$
where $Z_j^\tau = \alpha^\tau_{i_j} (X)$  when
$\tau = (\alpha^\tau_{i_1}, \dots, \alpha^\tau_{i_r} )$.
\end{enumerate}
\end{defn}

\begin{remark} \label{r:rembas}
Note that a diagonal basis is  not a basis for ${\bf V}$, but
rather a collection of bases for ${\bf V}$ (which give a basis for
the  distinguished subset $S_\Delta$ of the meromorphic functions on ${\bf V}$
determined by the hyperplanes in $\bigtriangleup$).
\end{remark}

\begin{remark} \label{r:remdiagord} (\cite{BV1},
Proposition 14; \cite{Sz}, Proposition 3.4)
A total order on $\tri$ gives rise to a diagonal
basis.
\end{remark}

\newcommand{\emm}{\mbox{\bf M}}
\newcommand{\enn}{\mbox{\bf N}}

\begin{theorem} \label{t:szgen}
[Szenes' theorem] \cite{Sz,BV2}
Let $f$ be a meromorphic function  on $\liet\otimes \CC$ with poles only on a
collection $\bigtriangleup$ of hyperplanes. Let $\enn$
be a lattice in $\liet$ and $\enn_{\rm reg}$ the regular points
in $\enn$. Let
$\emm$ be the lattice dual to $\enn$. Let $OB$ be a diagonal basis
associated to the hyperplane arrangement $\bigtriangleup$.
Then
\begin{equation} \label{e:1.2}
\sum_{n \in \enn_{\rm reg} } e^{<t, 2 \pi i n>} f(2 \pi in)
= \sum_{\sigma \in OB} {\rm Res}^{\sigma} \left (
f(z) F^{t }_\sigma(-z) \right ) .
\end{equation}
\end{theorem}
Here, for $\sigma \in { OB}  $ we define
\begin{equation}  \label{e:1.2b}
F^t_\sigma (-z) = \frac{1}{|\emm/\emm_\sigma| }
\sum_{m \in {\mathcal{R} } (t, \sigma) }
\frac{e^{(t-m)(z) } }{\prod_{\alpha \in \sigma} (1 - e^{\alpha(z)} )  }
\end{equation}
where ${\mathcal{R} } (t, \sigma) =
\{ u \in \emm |~ t -u = \sum_{\alpha \in \sigma} n_\alpha \alpha, ~
0 \le n_\alpha < 1 \}. $
It follows that
${\mathcal{R} } (t, \sigma) $ is in bijective correspondence
with $
\emm/\emm_\sigma$, where
$\emm_\sigma = \oplus_{\alpha \in \sigma} {\mathbb{Z}} \alpha
\subset \emm. $

\begin{remark}
Notice that the left hand side of (\ref{e:1.2}) is independent of the
choice of diagonal basis $OB$.
\end{remark}

\newcommand{\intlat}{\Lambda^I}

\begin{remark}
For our purposes $\enn$ is usually the integer lattice
$\intlat\subset \liet $ and $\emm$ the weight lattice.
\end{remark}

\subsubsection{$G=SU(n)$}
 The Lie algebra $ {\bf{t}} = {\bf{t}}_{n-1} $
 of the maximal torus $T$ of $SU(n)$ is
\begin{equation} \label{e:hypn}
   {\bf{t}} = \{(X_1,...,X_n) \in {\mathbb R}^n : X_1 + \cdots + X_n = 0 \}
\end{equation}
 Define coordinates $ Y_j = e_j(X) = X_j - X_{j+1} $ on $\bf{t}$ for $j=1,...,n-1.$
 The positive roots of $SU(n)$ are then $\gamma_{jk}(X) = X_j - X_k = Y_j + \cdots + Y_{k-1} $
 for $ 1 \le j < k \le n. $
 The {\em integer lattice}  $\Lambda^I$ of $SU(n)$ is generated by
 the simple roots $e_j, j=1,...,n-1.$ The dual lattice to $\Lambda^I$ with respect to the inner
 product $\langle \cdot, \cdot \rangle$
(in this case the Euclidean inner product on ${\mathbb R}^n$
restricted to the hyperplane (\ref{e:hypn}))
is the {\em weight lattice}
  $\Lambda^w \subset {\bf{t}}$
 \nolinebreak[4]; in terms of the inner product $\langle \cdot, \cdot \rangle$, it is given by $ \Lambda^w
 =\{ X \in {\bf{t}} : Y_j \in {\mathbb Z} \ for \ j = 1,...,n-1 \}. $ We define also $ \Lambda_{reg}^w
  = \{ X \in \Lambda^w :
 \ \gamma_{jk}(X) \not= 0
 \ {\rm for \ any} \ j \not= k \}. $

\newcommand{\diag}{{\rm diag}}

 \begin{defn}

  Let $ f : {\bf{t}} \otimes {\mathbb C} \rightarrow {\mathbb C} $
 be a meromorphic function$^1$\footnote{1. We shall use this notation
for a particular class of rational functions $g$.} of the form
  \[
    f(X) = g(X)e^{-\gamma(X)}
  \]
  where
$ \gamma(X) = \gamma_1Y_1 + \cdots + \gamma_{n-1}Y_{n-1} $ for $ (\gamma_1,...,\gamma_{n-1}) \in
  \mathbb{R}^{n-1}. $ We define
\begin{equation} \label{e:doublesqbrack}
    [[ \gamma ]] = ( [[\gamma]]_1,...,[[\gamma]]_{n-1})
\end{equation}
  to be the element of $ {\mathbb R}^{n-1} $ for which $ 0 \le [[\gamma]]_j < 1 $ for all $ j = 1,...,n-1 $
  and $ [[\gamma]] = \gamma$  $mod \ {{\mathbb Z}^{n-1}}$.
(In other words, $ [[\gamma]] = \sum_{j=1}^{n-1}
  [[\gamma]]_j e_j $ is the unique element of $ {\bf{t}}
\cong {\mathbb R}^{n-1} $ which is in the fundamental
  domain defined by the simple roots for the translation action on
$ {\bf{t}} $ of the integer lattice,
  and which is equivalent to
 $\gamma$ under translation by the integer lattice.)

  We also define the meromorphic function $ [[f]] : \bf{t} \otimes
{\mathbb{C}} \rightarrow {\mathbb{C}} $ by
  \[
    [[f]](X) = g(X) e^{-[[\gamma]](X)}.
  \]

 \end{defn}

 \begin{theorem} \label{t:szsun}
  [Szenes' theorem for $SU(n)$] (\cite{JK3},\cite{Sz})  \label{t:szenes}
  Let $ f : {\bf{t}} \otimes {\mathbb C} \rightarrow {\mathbb C} $
be a meromorphic
  function of the form $ f(X) = g(X) e^{- \gamma(X)} $ where $ \gamma (X) =
  \gamma_1 Y_1 +,...,+ \gamma_{n-1} Y_{n-1} $ with $ 0 \le \gamma_{n-1} < 1 $, and $g(X)$ is a rational
  function of $X$ decaying rapidly at infinity with poles only at
the zeros of the roots $\gamma_{jk}$. Then
\begin{equation} \label{e:t1.5}
    \sum_{\lambda \in \Lambda_{reg}^w ({\bf{t}})} f(2 \pi i \lambda) = Res_{Y_1=0}...Res_{Y_{n-1}=0}
    \left( \frac{ \sum_{w \in W_{n-1}} [[w (f)]] (X) } {(e^{-Y_{n-1}}-1) \cdots (e^{-Y_1}-1)} \right)
  \end{equation}
  where $W_{n-1}$ is the Weyl group of $SU(n-1)$ embedded in $SU(n)$ using the first $n-1$ coordinates.
\end{theorem}

 \begin{defn} \label{r:1.4}
Let $c$ be an element of $Z(G)$ (the center of $G$).
\end{defn}

\begin{remark} Note that Theorem \ref{t:szsun} refers to
a particular diagonal basis for $SU(n)$ which
does not generalize to other Lie groups.
See Section 5.3 of  \cite{Sz}. Applying Theorem \ref{t:szgen}
to the case $G = SU(n)$ we obtain an alternative formula for the
left hand side of (\ref{e:t1.5}).
\end{remark}

 \subsection{Residue formulas for the Hamiltonian space}

Let $N$ be a symplectic manifold with a Hamiltonian
$G$-action. We denote the moment map by
$\mu: N \to {\bf g}^*$, and define
the symplectic reduced space
\begin{equation} \label{e:reddef} N_{\rm red}= \mu^{-1}(0)/G. \end{equation}
(Note that $\mu$ will refer to moment maps for Hamiltonian
group actions, while $\Phi$ refers to moment maps for
quasi-Hamiltonian group actions.)
The symplectic form on $N_{\rm red}$ will be denoted
$\omega_{\rm red}$.
We define
\begin{equation} \label{e:barwdef}
{\bar{\omega}}(X) = \omega + \mu (X)
\end{equation}
Then ${\bar{\omega}} \in \Omega_G^2(N)$
 and it is closed under the
Cartan model differential  and thus defines an element $[{\bar{\omega}}]
 \in H_G^2(N)$. Here $\Omega^2_G(N)$ refers to the collection of elements
of degree 2 in the Cartan model, while $H^2_G(N)$ is the
equivariant cohomology of $N$ with complex coefficients.
See for example \cite{BGV}.

Suppose $0$ is a regular value of $\mu$. Then
there is a natural map $\kappa : H_G^*(N) \rightarrow H^*(N_{red})$ defined by
 \[
   \kappa : H_G^*(N) \to H_G^*( \mu^{-1}(0) ) \cong H^*(N_{red}).
 \]
 It is obviously a ring homomorphism. A similar map is defined
for quasi-Hamilto\\-nian $G$-spaces, which we shall also denote by
$\kappa$.

The residue theorem (\cite{JK1}, Theorem 8.1) gives a formula for the
evaluation of $\kappa(\eta)$ on the fundamental class of $N_{\rm red}$.
 Guillemin and Kalkman and independently Martin have given an alternative version of the
 residue formula which uses the one-variable result inductively:

 \begin{theorem} \label{t:gkm}
   (\cite{GK} Theorem 3.1; \cite{martin})
  Suppose $N$ is a compact symplectic manifold acted on by a torus $T$
 in a Hamiltonian fashion with moment map $\mu$, and $\eta \in H_T^*(N)$.
The corresponding Kirwan map is denoted by
$\kappa_T: H^*_T(N) \to H^*(N_{red})$. Assume $0$ is a regular value of $\mu$.
Then
$$
   \int_{N_{red}} \kappa_T(\eta) =  {\sum_{N_1}}' \int_{(N_1)_{red}} \kappa_{T/T_1}
(Res_1 \eta).
$$
 Here, $N_1$ is a component of  the fixed point set of a circle subgroup
$T_1\cong S^1$ of $T$ (so that $\mu(N_1)$
 are critical values of $\mu$):
it is a symplectic manifold equipped with a Hamiltonian action
 of $T/T_1$ and the natural map
$\kappa_{T/T_1} : H_{T/T_1}^*(N_1) \rightarrow H^*((N_1)_{red})$, where
 $(N_1)_{red}$ is the reduced space
$$ \frac{(\mu_{T/T_1})^{-1}(0)} { {(T/T_1)} }$$
for the action of $T/T_1$ on $N_1$.
  We choose a ray $\lambda$
from $0$ to the
complement of the moment map image of $N$, and
sum over those components
$N_1$ whose moment map image intersects $\lambda$.
 The map $Res_1 : H_T^*(N) \rightarrow H_{T/T_1}^*(N_1)$ is defined in
\cite{GK}
  (3.6) as
\begin{equation} \label{e:defres}
   {\rm Res}_1 (\eta) = Res_{X_1 = 0}(\frac{i_{N_1}^* \eta}{e_{N_1}})
\end{equation}
 where
$$
   i_{N_1}^* \eta \in H_T^*(N_1) = H_{T/T_1}^*(N_1) \otimes H_{T_1}^* $$
is induced from the inclusion map $ i_{N_1}: N_1 \to N$,
 $e_{N_1} $ is the equivariant Euler class of the normal bundle
to $N_1$ in $N$,
 and $X_1 \in {\bf t}_1$ is a basis element for ${\bf t}_1$.
We have introduced the notation $H^*_T$ to denote
$H^*_T({\rm pt})$.
 \end{theorem}

The notation $\sum'$ means the sum
 over those
$T_1$ and $N_1$ for which a (generic) ray $\lambda$ in ${\bf t}^*$
 from $0$ to the complement of $\mu_T(N)$ intersects $\mu_T(N_1)$.
We will use this theorem to allow us to make an inductive argument.

 \section{Construction}

Let $M$ be a quasi-Hamiltonian $G$-space, where $G $
is a  compact Lie group,
 with moment map $\Phi: M \to G$ and $\omega \in \Omega^2(M)^G$.
Recall that in
Definition \ref{r:1.4}  we had chosen an element $c \in Z(G)$.

We can construct a corresponding Hamiltonian $G$-space $\widetilde{M}$ as follows.
Let ${\widetilde{M}}
= \{ (m, \Lambda) \in M \times {\bf{g}} | \Phi(m) = c \exp (\Lambda) \} $
and $\mu: M \times {\bf g} \to {\bf g}$
is defined by
\[
    \mu(m,\Lambda) = -\Lambda.
\]
Then the following  diagram commutes:

 \begin{equation} \label{e:commdiag}
   \begin{CD}
   \widetilde{M}     @>{-\mu}>> {\bf{g}}  \\
   \pi_1@VVV             @VVV c\exp   \\
   M         @>{\Phi}>>  G
   \end{CD}
\end{equation}

Since $d \exp^*\chi = 0$, we can find  $\sigma \in \Omega^2({\bf{g}})$ such that $ d \sigma =  \exp^*\chi $.
We see that $d(\pi_1^* \omega - \mu^* \sigma) = 0$.
The space  $\widetilde{M}$ becomes a Hamiltonian $G$-space with moment
map $\mu$ and the invariant 2-form
\begin{equation} \label{e:twoform}
{\tilde{\omega}} = \pi_1^* \omega - \mu^* \sigma \in \Omega^2({\widetilde{M}})^G . \end{equation}

Let $\eta  \in H^*_G(M).$
We want to evaluate
\[
  \int_{\Phi^{-1}(c)/G} \kappa(\eta e^{\bar{\omega}})
\]
>From our construction, we see $\Phi^{-1}(c) = \mu^{-1}(0)$.
 Therefore
\[
  \int_{\Phi^{-1}(c)/G} \kappa(\eta e^{\bar{\omega}})
  = \int_{\mu^{-1}(0)/G} \kappa(\eta e^{\bar{\omega}})
\]
(We have assumed $\Phi$ is proper, so these spaces are compact.)
So we will evaluate the integral over $\mu^{-1}(0)/G$.

We face the difficulty that
 $\widetilde{M}$ may be singular.
In order to overcome this problem we will use the {\em Poincar\'e
dual} introduced in
the next section.

 \section{Equivariant Poincar\'e Dual}
Since we know that $M \times {\bf{g}}$ is always smooth, we will work with integration over
$M \times {\bf{g}}$, instead of working with integration over $\widetilde{M}$.

 \medskip

 \begin{lemma}  \label{l3.1} (\cite{JK3} Corollary 5.6)
  Let $T$ be the maximal torus of $G$ acting on $G$ by
conjugation. If $c \in T$ then we
  can find a $T$-equivariantly closed differential form
$\hat{\alpha} \in \Omega_T^*(G)$ on $G$ with
  support arbitrarily close to $c$ such that
  \[
    \int_G \eta \hat{\alpha} = \eta |_c \in H_T^*
  \]
  for all $T$-equivariantly closed differential forms $\eta \in \Omega_T^*(G)$.
 \end{lemma}

\begin{remark} The proof of Lemma \ref{l3.1}
given in \cite{JK3} is for the case
$G = SU(n)$. Similar arguments apply for general $G$.
\end{remark}

 \begin{proposition} \label{p3.2}  (\cite{JK3} Proposition 5.7)
 Let $P : M \times {\bf g} \rightarrow G$ be
defined by
 \begin{equation} \label{e:pdef}
P: (m,\Lambda) \mapsto \Phi(m) \exp(-\Lambda),
\end{equation}
and let $c \in T$, so that
 ${\widetilde{M}} = P^{-1}(c)$.
Let $\alpha = P^* \hat{\alpha} $
(where $\hat{\alpha}$ was defined in Lemma \ref{l3.1}).
Hence for $\eta \in \Omega^*_G(M \times \lieg)$,
 \[
   \int_{M \times {\bf{g}}} \eta \alpha = \int_{{\widetilde{M}}} \eta
 \]

\end{proposition}

\begin{defn} \label{r:quottdef}
Let $N$ be a Hamiltonian $G$-manifold with moment map $\mu: N \to
{\bf g}^*$. We  introduce   the notation $N \quott G$ to denote
the symplectic quotient $\mu^{-1}(0)/G$.$^2$\footnote{2. This
notation is borrowed from algebraic geometry, where it refers to
the geometric invariant quotient of a K\"ahler manifold by the
action of a reductive group, which is identified with the
symplectic quotient of the maximal compact subgroup. Strictly
speaking the notation does not apply to our situation since the
objects we work with are not usually K\"ahler manifolds.}
\end{defn}
We have
${\widetilde M}\quott G = (P^{-1}(c) \cap \mu^{-1}(0))/G$.

\newcommand{\openset}{{\mathcal V}}

  Let $V$ be a small neighbourhood of $c$ in $G$.
  In fact, if $V^{'}$ is any neighbourhood of $c$ in $T$ containing
the closure of ${ V}$ then$^3$\footnote{3. After Lemma \ref{l3.3b}
we shall replace $V$ by a space  $\openset$ constructed there with
$c \in \openset$.}
  \[
    \int_{{P^{-1}(V^{'})}} \eta \alpha = \int_{\widetilde M} \eta \in H_T^*
  \]

 \begin{prop}
  {\bf{(Reduction to the abelian case)}} [S. Martin \cite{martin}]
\label{pr:redab}
  If $T$ is a maximal torus of $G$ and $M$
is a symplectic manifold equipped with
an effective Hamiltonian action of $G$, then denote by $\mu: M \to \lieg$
the moment map for $G$ and $\mu_T: M \to T$ the $T$ moment map.
Assume $0$ is a regular value of $\mu$. Then  for any regular value
  $\xi$ of $\mu_T$ sufficiently close to $0$ and $\eta \in
H^*_T(M)^W$  we have that
   \begin{equation} \label{6.1}
    \int_{\mu^{-1}(0)/G} \kappa^0(\eta  e^{\bar{\omega}} )
         = \frac{n_1}{n_0^{'}|W|} \int_{\mu^{-1}(0)/T}
\kappa^0_T ({\mathcal D} \eta e^{\bar{\omega} })  \end{equation}
\begin{equation} \label{6.2}
         = \frac{n_1}{n_0^{'}|W|} \int_{\mu^{-1}(\xi)/T}
\kappa^\xi _T ({\mathcal D}\eta e^{\bar{\omega} }) \end{equation}
 \begin{equation} \label{6.3}
        = \frac{(-1)^{n_+}n_1}{n_0^T |W|} \int_{\mu_T^{-1}(\xi)/T}
\kappa_T^\xi({\mathcal D}^2 \eta e^{\bar{\omega}})  \end{equation}
  where $n_1$ is the order of the stabilizer in $G$ of a generic point of $\mu^{-1}(0)$ and $n_0^T$
  (respectively $n_0^{'}$) is the order of the stabilizer in $T$ of a generic point of $\mu_T^{-1}(0)$
  (respectively $\mu^{-1}(0)$).
Here
 the polynomial ${\mathcal D} : {\bf t} \rightarrow {\mathbb R}$
is defined by
\begin{equation} \label{e:ddef}
{\mathcal D} (X) = \prod_{\gamma > 0} \gamma (X) \end{equation}
The equalities (\ref{6.1}) and (\ref{6.2}) are valid when
the restriction of $\eta$ to $\mu^{-1} (\liet)$ has compact support, or  when
$\mu$ is proper.
The equality (\ref{6.3}) is valid provided
 $\mu_T$ is proper.
\end{prop}

\begin{remark} The equality (\ref{6.3}) is proved in
\cite{JK1}, where explicit formulas for the left and right hand
sides are proved (establishing the equality).
\end{remark}

\begin{remark} \label{r:slice}
If a symplectic manifold $N$ is acted on
by a compact Lie group $G$ with maximal
torus $T$, we use the fact that
$$\mu^{-1}(0) = \mu^{-1}(\liet) \cap \mu_T^{-1}(0)$$
to apply Theorem \ref{t:gkm} when $G$ is nonabelian.
Note that $\mu^{-1}(\liet)$ is smooth because it can be written
as $(\mu^\perp) ^{-1} (0)$ where $\mu^\perp = \pi^\perp \circ \mu$
for $\pi^\perp: \lieg \to \liet^\perp$ the projection,
and we can easily
prove $d \mu^\perp$ is surjective because $\mu$ is equivariant.
Although the hypothesis  that
$\mu: N \to \lieg $ is proper does not
imply that $\mu_T: N \to \liet$ is proper,
it  does imply that
$\mu_T: \mu^{-1}(\liet) \to \liet$ is proper.
\end{remark}

\begin{remark}
By using Proposition \ref{pr:redab} we can relax the condition
from Definition \ref{r:1.4} that
$c \in Z(G)$ to $c \in T$ satisfying a constraint which implies that
$\Phi^{-1}(c)/T$ is an orbifold. See  Lemma \ref{l3.3b} below.
\end{remark}

 \begin{lemma} \label{l3.3a} (cf. \cite{JK3}, Lemma 5.10)
Let $c \in T$ be a regular value of $P$.
Then $(P^{-1}(V) \cap \mu^{-1}(0))/T$ is
  an orbifold (for $V$ as above).
\end{lemma}

\proof  Our first observation is that near
$(\Phi\circ \pi_1)^{-1}(c)$, $M$ is endowed with a  symplectic
structure, since
there is a $G$-invariant neighbourhood $V \subset G$
containing $c$  such that
the restriction of the closed 2-form
$\tilde{\omega}$ (defined in  (\ref{e:twoform})) to
$(\Phi\circ \pi_1)^{-1}(V)$ is nondegenerate.
This is true for the following reason.
Consider the diagram (\ref{e:commdiag}).
The space $M$ is a smooth manifold,
so the space $\widetilde{M}$ is smooth
on the preimage under $\Phi \circ \pi_1$ of all
regular values of
$c\exp$
(these regular values contain a
neighbourhood $V\subset G$ of $c$).
The two-form $\tilde{\omega}$
is closed on
$\widetilde{M}$.
The map $\mu$ satisfies the moment map condition
$ d\mu_\xi  = \tilde{\omega} (\nu_\xi, \cdot) $ on
$\widetilde{M}$.
(See \cite{AMM} and \cite{J}.)
Furthermore the 2-form $\tilde{\omega}$ descends
under symplectic reduction from $(\Phi \circ \pi_1)^{-1}(V)$
 to the standard symplectic form on
$\mu^{-1}(0)/G$.
It follows that $\tilde{\omega}$ is nondegenerate
on $(\Phi\circ \pi_1)^{-1}(V), $
which is an open neighbourhood of $\mu^{-1}(0)$ in
$\widetilde{M}$.

We know that $c$ is a regular value for
$P : M \times {\bf g} \to G$,
and therefore we can choose the neighbourhood
$ V$  so that all points of $V$ are also regular
values of $P$
(by standard properties of  the rank of a differentiable map).
Because $M \cap \Phi^{-1}(V)$ is symplectic,
with moment map $\mu$ related to $\Phi$ as in (\ref{e:commdiag}),
$G$ acts with finite stabilizers at all points of
$ (\Phi\circ \pi_1)^{-1}(V) \cap \mu^{-1}(0)$.
This implies that $T$ also acts with finite stabilizers at
all points of $ (\Phi\circ \pi_1)^{-1}(V) \cap \mu^{-1}(0)$.
Hence $P^{-1}(V) \cap \mu^{-1}(0)/T$ is an orbifold.
\hfill $\Box$

 We define a one dimensional torus ${\hat{T}}_1 \cong S^1$ in $G$
generated by a chosen  element ${\hat{e}}_1 $ in the integer lattice of the
Lie algebra of a fixed maximal torus
$T$. The corresponding coordinate will be denoted
\begin{equation} \label{e:yonedef}
Y_1 = \langle \hat{e_1}, X \rangle.
\end{equation}
Then $\hat{T_1}$
 is identified with $S^1$ via
 \[
   e^{2 \pi i t}  \in S^1 \mapsto \exp{t \hat{e}_1 }
\in {\hat{T}}_1
 \]
 The one dimensional Lie algebra
${\hat{\bf{t}}}_1$ is spanned by ${\hat{e}}_1$. Its
 orthocomplement in ${\bf t}$ (under the inner product)
is defined as ${\bf t}_{r-1}$.
 Define ${\widetilde T}_{r-1}$
to be the torus given by $\exp ( {\bf t}_{r-1} )$.
Also define $T_{r-1}$ to be the quotient torus
  $T_{r-1} =  T/\hat{T_1},$ so that its Lie algebra is  also
${\bf t}_{r-1}$.

 \begin{lemma} \label{l3.3b} (cf. \cite{JK3},  Proposition 6.5)
There is a $T$-invariant open neighbourhood $\openset \subset G$ for which
the group  $T_{r-1}$ acts locally freely on
$P^{-1}({\openset}) \cap \mu^{-1}({\hat{\liet_1}})$,
so $P^{-1}(\openset) \cap \mu^{-1}(\hat{\liet_1})/T_{r-1}$ is an orbifold.
\end{lemma}

\proof
Consider the map $P: M \times \lieg \to G$ defined at
Proposition \ref{p3.2}.
This is a $G$-equivariant map (where $G$ acts on $\lieg$ and
$G$ by the
adjoint action).

\newcommand{\subsym}{{\mathcal W}}
\newcommand{\locfr}{{\mathcal U}_{r-1}}
\newcommand{\locg}{{\mathcal U}_G}
Let $\locfr \subset M$ be the subset of $M$ where
$T_{r-1}$ acts locally freely.
\newcommand{\subx}{{\mathcal W}_x}
There is an open dense subset $\subsym$ of $M$  for which
open neighbourhoods $\subx$ of $x$ are equipped with
 Hamiltonian $G$ actions
$ \mu: \subx \to \lieg$. (See \cite{AMM}.)

Let
$$M_{G} = \{ m \in M ~|~ G
~\mbox{acts locally freely at~}m\}.$$
This set is $T$ invariant.

We note that $d\Phi: T_x M_G \to T_{\Phi(x)} G$ is
surjective for all $x \in M_G$, and  by the inverse function theorem
this implies $\Phi(M_G)$ is an open  subset of $G$.

Hence $P (M_G \times \liet_1) $ is also
an open  subset of $G$.

We define $ \openset =  P(M_{G} \times \liet_1) . $
This is a $T$-invariant open
set.
>From now on, without loss of generality we replace our
earlier open neighbourhood $V \subset G$
(introduced after Definition \ref{r:quottdef})
by $\openset$, and
we redefine $c$ to be a point in $\openset $.

\hfill $\Box $

We extend the definition of the composition
  \[
    \kappa : H_T^*(P^{-1}(c)) \rightarrow H_T^*(P^{-1}(c) \cap \mu^{-1}(0))
           \cong H^*(P^{-1}(c) \cap \mu^{-1}(0)/T)
  \]
  to
  \[
    \kappa :  H_T^*(P^{-1}(\openset)) \rightarrow H_T^*(P^{-1}(\openset)
 \cap \mu^{-1}(0))
           \cong H^*(P^{-1}(\openset) \cap \mu^{-1}(0))/T).
  \]
By Proposition \ref{p3.2} we have

 \[
  \int_{P^{-1}(c) \cap \mu^{-1}(0)/T} \kappa(\eta)
         = \int_{P^{-1}(\openset) \cap \mu^{-1}(0)/T} \kappa(\eta \alpha)
 \]
where $\eta \in H^*_G(M)$.
 The class $\alpha$ is the  Poincar\'e  dual
of $P^{-1}(c)$ in $P^{-1}(\openset)$.

 \section{Periodicity}

The space $\widetilde{M} $ is  noncompact, so many standard results
in symplectic geometry cannot be applied without further  analysis.
To handle the difficulties of noncompactness
we rely on the fact  that the moment map $\mu$ is proper
(because $\Phi$ is proper and
$\mu^{-1}(U) = \Phi^{-1}(\exp U)$) and use a relative version of
Guillemin-Kalkman's theorem (Proposition \ref{p:4.2}) which is valid
provided the inverse image of a closed interval under the moment map
for a circle action
is compact. Using Proposition \ref{p:4.2} we reduce to
the components of the fixed point set of the circle action
 for which the moment map
takes values in a closed interval. Since the moment map for this
circle action
is proper there are only finitely many components of the fixed point
set whose moment map takes values in this closed interval, so
the sums in
Proposition \ref{p:4.2} and Lemma \ref{l:4.1} are finite sums.

 \begin{prop}  \label{p4.1} (\cite{JK3}, Proposition 6.3;\cite{Goldin})
  For any symplectic manifold $M$ equipped with a
Hamiltonian action of $T = T_r$ such that ${\widetilde T}_{r-1}$
  acts locally freely on $\mu_{{\widetilde T}_{r-1}}^{-1}(0)$,
the symplectic quotient $\mu_{T_r}^{-1}(0)/T$
  may be identified with the symplectic quotient of
$\mu_{{\widetilde T}_{r-1}}^{-1}(0)/{\widetilde T}_{r-1}$ by the
  induced Hamiltonian action of ${\hat{T_1}}$. Moreover if in addition
$T_r$ acts locally freely
  on $\mu_{T_r}^{-1}(0)$ then the ring homomorphism
$\kappa : H_{T_r}^*(M) \rightarrow
  H^*(\mu_{T_r}^{-1}(0)/T_r)$ factors as
 \begin{equation} \label{e:factor}
    \kappa = {\hat{\kappa}_1} \circ \kappa_{r-1}
\end{equation}
  where
  \[
    \kappa_{r-1} : H_{T_r}^*(M) \rightarrow
H_{T_r}^*(\mu_{{\widetilde T}_{r-1}}^{-1}(0))
             \cong H_{{\hat{T}}_1 \times {\widetilde T}_{r-1}}^*
(\mu_{{\widetilde T}_{r-1}}^{-1}(0))
             \cong H_{{\hat{T}}_1}^*(\mu_{{\widetilde T}_{r-1}}^{-1}(0)
/{\widetilde T}_{r-1})
  \]
  and
  \[
   {\hat{\kappa}_1} :  H_{{\hat{T}}_1}^*(\mu_{{\widetilde T}_{r-1}}^{-1}(0)
/{\widetilde T}_{r-1})
\rightarrow
                           H^* \Bigl ( (\mu_{{\widetilde T}_{r-1}}^{-1}(0) \cap
\mu_{{\hat{T_1}}}^{-1}(0))/({\widetilde T}_{r-1} \times {\hat{T}}_1 ) \Bigr )
                       \cong H^*(\mu_{T_r}^{-1}(0)/T_r)
  \]
are the corresponding compositions of restriction maps with
similar isomorphisms (i.e., Kirwan maps).
 \end{prop}

\begin{remark} \label{r4.2}
Let $T_r$ and $T_{r-1}$ be as above Proposition \ref{p4.1}.
Note  also that
\begin{equation} \label{e:factor2}
\kappa = \hat{\kappa}_{r-1} \circ {\kappa_1} \end{equation}
where
$${\kappa_1}: H^*_T(M) \to H^*_{(T/\hat{T_1})} (\mu_{\hat{T_1}}^{-1}(0)/\hat{T_1})$$
and
$$\hat{\kappa}_{r -1}: H^*_{(T/\hat{T_1})} (\mu_{\hat{T_1}}^{-1}(0)/\hat{T_1} )
\to H^*_T(\mu_T^{-1}(0))\cong H^*(\mu_T^{-1}(0)/T)$$
are defined in a way similar to the definitions of
$\hat{\kappa}_1$ and $\kappa_{r-1} $ given in the
statement of the previous proposition.
\end{remark}

We introduce the notation $\kappa^0: H^*_G(M) \to
H^*(\mu^{-1}(0)/G) $ to mean  the composition
of the restriction map $H^*_G(M) \to H^*_G(\mu^{-1}(0))$ and
the isomorphism
$H^*_G(\mu^{-1}(0)) \cong H^*(\mu^{-1}(0)/G)$.
The map
$\kappa^0_T: H^*_T(M) \to H^*(\mu_T^{-1}(0)/T)$ is
defined similarly, and
in the same way
$\kappa^0_T:  H^*_T(M ) \to H^*_T(M \cap \mu^{-1}(\liet))
\to H^*(M \cap \mu^{-1}(0)/T)$ (since
$\mu^{-1}(0) = \mu_T^{-1}(0) \cap \mu^{-1}(\liet)).$
The maps $\kappa^\xi_T$ are defined in  the same way,
except that $\mu^{-1}(0)$ is replaced by $\mu^{-1}(\xi)$.

 \begin{prop}
\label{p:4.2}
 {\bf{(Dependence of symplectic quotients on parameters) [Guillemin-Kalkman \cite{GK}; S. Martin \cite{martin}] }}
Let $M$ be a  Hamiltonian $T$-space for which the moment map $\mu_T$ is
proper.
 If $T = U(1)$ and $n_0^T$ is the order of the stabilizer in $T$ of a generic point of $\mu_T^{-1}(0)$ then
\begin{align*}
   \int_{\mu_T^{-1}(\xi_1)/T} \kappa^{\xi_1}_T (\eta e^{\bar{\omega}})
   - \int_{\mu_T^{-1}(\xi_0)/T} \kappa^{\xi_0}_T(\eta e^{\bar{\omega}})
  =&\\& \hskip -60pt= n_0^T \sum_{E \in {\mathcal E} : \xi_0 < \mu_T(E) < \xi_1} Res_{X=0} e^{\mu_T(E)X} \int_E \frac{\eta(X) e^\omega }{ e_E(X) }
 \end{align*}
 where $X \in {\mathbb C}$ has been identified with $2 \pi i X \in {\bf t} \otimes {\mathbb C}$
 and $\xi_0 < \xi_1$ are two regular values of the moment map.
Here ${\mathcal E}$ is the set of all components
$E$ of the fixed point set of $T$ on $M$.
 \end{prop}

 \begin{lemma} \label{l:4.1}  {( \bf cf. \cite{JK3} Lemma 6.7)} \label{period}
Let $M$ be a quasi-Hamiltonian $G$-space for which $\Phi$ is
proper,$^4$\footnote{4. The hypothesis that $\Phi$ is proper is
satisfied in many examples, including all compact examples.} and
let the space $\widetilde{M}$ be defined by (\ref{e:commdiag}).
Suppose $\eta \in H^*_G(M).$
   If $\openset$ is as defined in the proof of Lemma \ref{l3.3b},
and we define
$M = P^{-1}(c)$    for
some $c$ in $\openset$ so
that (by Lemma \ref{l3.3b})
$P^{-1}(\openset) \bigcap \mu^{-1}({\hat{\bf t}}_1)/T_{r-1}$
 is an orbifold and we define
  $N(\openset) = P^{-1}(\openset) \bigcap \mu^{-1}(0)/T$, then
this hypothesis implies also that  $0$ is a regular value of
$\mu_{T/\hat{T_1}}: E \to \liet/\hat{\liet_1} $ for all $E \in {\mathcal E}$,
  where ${\mathcal E}$ is the set of
components of the fixed point set of the action of ${\hat{T}}_1$ on
$M \times \hat{\liet_1}$.
We then have that
\begin{equation} \label{e:6.7statement}
    \int_{N(\openset)} \kappa(\eta e^{\bar{\omega}} e^{-Y_1} \alpha)
        = \int_{ P^{-1}(\openset) \bigcap \mu^{-1}({\hat{e}}_1)/T}
\kappa(\eta e^{\bar{\omega}} \alpha)
\end{equation}
\begin{equation} \label{e:6.7statementbis}
=\int_{N(\openset) }
 \kappa (\eta e^{\bar{\omega} } \alpha)
    -n_0
\sum_{ E \in {\mathcal E} :
-|| {\hat{e}_1 }||^2 < ~\langle {\hat{e}_1},\mu(E) \rangle <0  }
           \int_{E\quott \widetilde{T}_{r-1} }
   \kappa_{T/\hat{T_1} }  Res_{Y_1 = 0} \frac{
\eta e^{\bar{\omega}} \alpha } {e_E}.
\end{equation}
Here
$e_E$ is the $\hat{T_1}$-equivariant Euler
class of the normal bundle
to $E$ in ${\widetilde{M}}$,
 while $n_0$
is the order of the subgroup of ${\hat{T}}_1/({\hat{T}}_1 \bigcap
{\widetilde{T}_{r-1}})$
  that acts trivially on
$P^{-1}(\openset) \cap \mu^{-1}({\bf{\hat{t_1} } } )$, and $Y_1$ was defined in
(\ref{e:yonedef}).
 Also $\alpha$ is
  the $T$-equivariantly
closed differential form on $M \times {\bf t}$ given by Proposition \ref{p3.2}
  which represents the equivariant Poincar\'e  dual  of
${\widetilde M}$, chosen so that the
  support of $\alpha$ is contained in $P^{-1}(\openset)$.
\end{lemma}

\proof
Note that throughout this proof we make
extensive use of Remark \ref{r:slice}.
  Since $\mu^{-1}({\bf \hat{t}}_{1}) = M \times {\bf \hat{t}}_{1}$ is contained in
$\mu_{\widetilde{T}_{r-1} }^{-1} (0)$, it follows
  from Lemma \ref{l3.3b}
 that
 $P^{-1}(\openset) \cap \mu^{-1}
( \hat{\bf t}_1)/{\widetilde{T}_{r-1}}$ is an
orbifold.
(Here we let  $\mu_{\hat{T}_1} = \pi_{\hat{\bf t}_1} \circ \mu$, and
  $\pi_{\hat{\bf t}_1} : {\bf t} \rightarrow \hat{\bf t}_1$
is the  projection corresponding to the bi-invariant inner product.)
The restriction of $\mu_{\hat{T}_1}$ to $\mu^{-1}({\bf \hat{t}_1})$
is proper, because if $C \subset
\hat{\bf t}_1 $ is compact,
then $$\mu_{\hat{T}_1}^{-1}(C) =
\{ (x ,\Lambda) \in
M \times C: \Phi(x) =\exp(\Lambda) \}. $$
This is compact because $C$ is covered by finitely many compact
subsets $U$ of ${\hat{\liet}_1}$  on each of
which the restriction of the exponential
map to ${\hat{\liet}_1}$ is a diffeomorphism, and on each
 $U$ we find that $\{(x, \Lambda) \in M \times U:
\Phi(x) = \exp(\Lambda) \} = \Phi^{-1} (\exp U)$,
which is compact  because $\Phi$ is proper.
Also,
the support of $\alpha$ is contained in $P^{-1}(\openset)$.
To apply Proposition \ref{p:4.2}
we need that $0$ and $\hat{e_1}$ are regular values of
$\mu_T: \mu^{-1}(\liet) \to \liet$.
This follows from Lemma \ref{l3.3a} where it is shown that
$N(\openset)$ is an orbifold. (Note that by definition
$N(\openset)$ is identical to $P^{-1}(\openset) \cap \mu^{-1}(\hat{e_1})/T_r$.)

Hence
  Guillemin and Kalkman's proof of Proposition \ref{p:4.2}
 can be applied to the $\hat{T}_1$-invariant function induced by $\mu_{\hat
  {T}_1}$ on $\widetilde{M} \cap \mu^{-1}(\liet)$ and the $\hat{T}_1$-equivariant form induced by
  $\eta e^{\bar{\omega}} \alpha$. By combining this with
Remark \ref{r4.2} and using Theorem \ref{t:gkm},
Remark \ref{r:slice}  and (\ref{6.1}) of Proposition \ref{pr:redab} we get
$$
       \int_{N(\openset)} \kappa(\eta e^{\bar{\omega}} \alpha)
                - \int_{ P^{-1}(\openset)^T \bigcap \mu^{-1}({\hat{e}}_1)/T}
\kappa(\eta e^{\bar{\omega}} \alpha)
$$
\begin{equation} \label{e:6.7pf}
      = n_0 \sum_{ E \in {\mathcal E} : -||{\hat{e}}_1||^2 <
\langle {\hat{e}}_1,\mu(E) \rangle <0 }
              \int_{E \quott T_{r-1}}
  Res_{Y_1 = 0}\kappa_{r-1}
\frac{ (\eta e^{\bar{\omega} }  \alpha)}{e_E }
\end{equation}

 Now we need to show (\ref{e:6.7statement}).
 Note that the restriction of $P : M \times {\bf g} \rightarrow G$ to $\mu^{-1}({\bf t}) = M \times {\bf t}$ is invariant
 under the translation $s_{\Lambda_0} : M \times {\bf g} \rightarrow M \times {\bf g}$ defined by
 \[
     s_{\Lambda_0} : (m, \Lambda) \mapsto (m, \Lambda + \Lambda_0)
 \]
 for $\Lambda_0 \in \Lambda^I = \ker(\exp)$ in ${\bf t}$. So
for $\Lambda_0 = \hat{e_1}$
 \begin{eqnarray}
    \int_{ P^{-1}(\openset)^T \bigcap \mu^{-1}({\hat{e}}_1)/T_r} \kappa(\eta e^{\bar{\omega}} \alpha)
        & = & \int_{N(\openset)} \kappa(s^*_{\hat{e}_1}(\eta e^{\bar{\omega}} \alpha)) \nonumber \\
        & = & \int_{N(\openset)}
\kappa(\eta e^{\bar{\omega}} e^{-Y_1} \alpha) \nonumber
 \end{eqnarray}
This proves (\ref{e:6.7statement}).

It is proved in \cite{GK} (Proof of (3.10)) that
\begin{equation} \label{e:3.10pf}
\kappa_{T/\hat{T_1}}  Res_1 =
Res_1  \kappa_{T/{\hat{T_1}}}  \end{equation}
where  the notation
$Res_1 $ was introduced in (\ref{e:defres})
and  on the left hand side of (\ref{e:3.10pf}) we have
$$Res_1: H^*_T(\mu^{-1}(\liet) \cap \widetilde{M})
\to H^*_{T/\hat{T_1}} (\widetilde{M}_1\cap \mu^{-1}(\liet) )$$
(where $\widetilde{M}_1$ is the fixed point set of the action of
$\hat{T}_1$ on $\widetilde{M}$)
while on the right hand side we have
$$ Res_1: H^*_{T_1} \Bigl (\mu^{-1}(\liet) \quott (T/T_1)\Bigr )  \to
H^*(\mu^{-1}(\liet)\quott T). $$
This enables us to deduce (\ref{e:6.7statementbis}) from
(\ref{e:6.7pf}). \hfill $\Box$

\begin{remark}
The proof of Lemma 6.7 in \cite{JK3} was for the case $M = G^{2h},$
but the only property of this space  used is that $\Phi$ is proper
and the explicit identification of a suitable
set $\openset$ in \cite{JK3}, Proposition 6.5.
\end{remark}

We now come to a  lemma  which asserts that the fixed point set
of the action of a circle subgroup on a quasi-Hamiltonian
$G$-space is also a quasi-Hamiltonian space.
This result is a special case of \cite{AMW3},
Proposition 4.4. The result is valid for subgroups $S$ of
rank higher than $1$, but we state it only for $S \cong S^1$ since
this is the only case we need to make the induction work.
We include a proof
for completeness.

 \begin{lemma} \label{l:4.2}
    Let $(M, \omega, \mu)$ be a quasi-Hamiltonian $G$-space and
let  $
\hat{T_1} = S \cong S^1$ be a subgroup of $G$. If $H$ is the fixed point set of the
    adjoint action of $S$ on $G$ (the subgroup of $G$ consisting
of all elements which commute with  all elements of
$S$), then
$H$ is a Lie subgroup of $G$ and
$M^S$ is a quasi-Hamiltonian $H$-space.
\end{lemma}

\proof{    Since $M$ is a quasi-Hamiltonian $G$-space, it satisfies
the  three axioms given in Definition 1.1.
We need to check that
these axioms are also valid for the  $H$-space $M^S$. \\
   (1)  Let $\iota_S : M^S \rightarrow M$, and $\Phi_H = \Phi_G|_{M^S}.$
(Notice that on $M^S$ $\Phi_G$ takes values in $H$, by equivariance
of $\Phi_G$.)
Then
     \begin{eqnarray}
            d \iota_S^* \omega &=& \iota_S^*d \omega              \nonumber \\
                                               &=& -i_S^* \Phi_G^* \chi_G       \nonumber \\
                                               &=& -\Phi_H^* \chi_H                 \nonumber
    \end{eqnarray}
   (2) If $\xi \in {\bf h}$
    \begin{eqnarray}
        \iota (\nu _\xi) \iota_S^*\omega
                                & = & \frac {1}{2} \iota_S^* \Phi_G^* \langle\theta_G + \bar{\theta}_G, \xi \rangle   \nonumber \\
                                & = & \frac{1}{2} \Phi_H ^* \langle \theta_H + \bar{\theta}_H, \xi \rangle     \nonumber
    \end{eqnarray}
    \noindent
   (3) At each $x \in M^S$ and for $\xi \in {\bf h}$, we need to show that the kernel of $\iota_S^*\omega _x$  is given by
    \begin{equation}
     \ker \iota_S^*\omega _x  =  \{ \nu _{\xi}(x), \xi \in \ker(Ad_{\Phi_H(x)} + 1) \}.
    \end{equation}
\noindent We know that for $x \in M^S$, ${\rm Ker} \omega_x =
\{\nu_\xi, \xi \in \lieg, {\rm Ad}_{\Phi(x)} \xi = - \xi \}. $ If
$v \in T_x M^S$ satisfies $\omega(v,Y) = 0 $  for all $Y \in T_x
M$ then $v = \nu_\xi $ for some $\xi \in \lieg$ with ${\rm
Ad}_{\Phi(x)} \xi =\nolinebreak[4] - \xi. $ It suffices to prove
that $\nu_\xi(x) \in T_x M^S$ if and only if $ \xi \in {\bf  h}$.
It is obvious that if $ \xi \in {\bf h}$ then $\nu_\xi(x) \in T_x
M^S$. To prove the converse, we see that if $\xi \notin {\bf h}$
then (restricting to $\xi $ which are
 orthogonal to the Lie algebra of $  {\rm Stab}(x)$ in
the bi-invariant inner product, so that the map
$\xi \mapsto \nu_\xi(x)$ becomes bijective)
${\rm Ad}_s \xi \ne \xi$ for $s$ a generator of $S$.
Hence ${\rm Ad}_s \nu_\xi \ne \nu_\xi$ which
implies $\nu_\xi \notin T_x M^S$, as we wished to prove.

\hfill $\Box$
}

\begin{remark} If $T$ is a torus which is a
subgroup of $G$ and $M$ a quasi-Hamiltonian $G$-space,
then $\Phi^{-1}(T)$ is a quasi-Hamiltonian $T$-space.
Thus also $\mu^{-1}(\liet)$ is a Hamiltonian $T$-space.
This is proved  in \cite{AMM}.
\end{remark}

\begin{remark}
Under the hypotheses of Lemma \ref{l:4.2}, $T \subset H$ and
$\hat{T_1} \subset Z(H)$. Thus $T_{r-1} = T/\hat{T_1}$
is a group of rank $r - 1$ with a quasi-Hamiltonian  action
 on $M^{\hat{T_1}}$. This enables us to perform an inductive argument.
\end{remark}
 \bigskip

 From Lemma \ref{l:4.1}, we get
 \begin{eqnarray}
   \int_{N(\openset)} \kappa(\eta e^{\bar{\omega}} \alpha)
    - \int_{N(\openset)} \kappa(\eta e^{\bar{\omega}} e^{Y_1} \alpha)
    &=& \int_{N(\openset)} \kappa(\eta e^{\bar{\omega}} (1-e^{Y_1}) \alpha) \nonumber \\
    &\hskip -80pt=& \hskip -40pt n_0 \sum_{ E \in {\mathcal E} : -||{\hat{e}}_1||^2 < \langle {\hat{e}}_1,\mu(E) \rangle <0 }
             \int_{E\quott T_{r-1}} \kappa_{r-1}
 Res_{Y_1 = 0}\frac{ \eta e^{ \bar {\omega } } \alpha }{ { e_E }}
 \end{eqnarray}
Thus
\begin{proposition} \label{p:perone}
Under the hypotheses and in the notation
of Lemma \ref{l:4.1}, we have
 \begin{eqnarray} \label{N_n}
 \hskip 20pt  \int_{N(\openset)} \kappa(\eta e^{\bar{\omega}} \alpha)
   =  n_0 \sum_{ E \in {\mathcal E} :
 -||{\hat{e}}_1||^2 < \langle {\hat{e}}_1,\mu(E) \rangle <0 }
          \int_{E\quott T_{r-1}}
  \kappa_{r-1} Res_{Y_1 = 0}
\frac{ ( \eta e^{\bar{\omega} } \alpha) }{e_E(1-e^{Y_1} ) }
 \end{eqnarray}
\end{proposition}

 \section{Residue formula for quasi-Hamiltonian spaces}

\newcommand{\ZZ}{{\mathbb Z}}

We shall now list some properties of diagonal
bases which will enable us to use the general
form of Szenes' theorem (Theorem \ref{t:szgen})
  to relate our results with  Alekseev-Meinrenken-Woodward's
results in Section 6.

Let $\hat{\liet_1}$ be the Lie algebra of a circle subgroup
$\hat{T_1} \cong S^1$.
Let $\hat{e_1}$ be an element of the
weight lattice of $T$, which projects to
 a generator of $\hat{\liet_1}^*$.
We will study diagonal bases for $\liet\otimes \CC$; these will correspond to
a choice of weights. We refer to Definition
\ref{d:diagbas} and Remarks \ref{r:rembas} and \ref{r:remdiagord}.

Let $S$ be a subgroup of $T$.
The purpose of the following lemma is to show that diagonal bases for hyperplane
arrangements for the Lie algebra of
$S/\hat{T}_1$ may be converted into diagonal bases for a hyperplane
arrangement in  the Lie algebra of $S$ by adjoining a generator
of $\hat{\liet_1}$ to each ordered basis of ${\rm Lie} (S/\hat{T}_1)$.

\newcommand{\lies}{{ \mathfrak{s} }}

\begin{lemma} \label{l6.2old} Let $\{\alpha'_1, \dots, \alpha'_{N_1} \} $ be
a collection of elements in $\lies^*$, where $S \triangleleft T$ is
a subgroup of $T$ such that $\hat{T_1} \subset S$, and $\lies$ is the
Lie algebra of $S$.
Suppose $\sigma' = \{ \sigma'_1, \dots, \sigma'_{N_2}\} $
is a diagonal basis for a hyperplane arrangement
$\bigtriangleup' $ in $\lies/\hat{\liet_1}$ corresponding to a chosen
 total
ordering on  the hyperplane arrangement $\Delta'$
which comprises
the
weights for the action of $S/\hat{T_1}$  on all  tangent spaces
$T_F M_1$ for all components $F$ of the fixed point
set of the action of   $S/\hat{T_1}$ on $M_1$ for
some subgroup $S$ of $T$ satisfying
$\hat{T_1} \subset S$,  where
 $M_1$ is a component of   $ M^{\hat{T_1}}$.
 So each
$\sigma'_j$ is a basis for $\lies/\hat{\liet_1}$.
  Then
$\sigma = \{ \sigma_1' \cup \hat{e_1}, \dots,\sigma'_{N_2}\cup \hat{e_1} \} $
is a diagonal basis for the hyperplane arrangement
$\bigtriangleup' \cup \{ \hat{e_1}\} $ in $\lies$, which corresponds to
a subset of the weights of the action of $S$ on all tangent spaces
 $T_F M$ at all
 components $F$ of the fixed point set
of
the action of $S$ on $M$.
\end{lemma}

\proof
Suppose
$\sigma'_j = (\alpha'_{j_1},
\dots, \alpha'_{j_{\ell-1 }}) $ and
$\sigma'_k = (\alpha'_{k_1},
\dots, \alpha'_{k_{\ell-1 }}) $
are ordered bases which are members of a
 diagonal basis for a hyperplane arrangement in
$\lies/\hat{\liet_1}$.
(Here, $\ell$ is the rank of $S$, and $k = (k_1, \dots, k_{\ell-1})$ and
$j = (j_1, \dots, j_{\ell-1})$ are multi-indices.)
 We note that if $\sigma_j' = (\alpha'_{j_1},
\dots, \alpha'_{j_{\ell-1 }}) $ and
$(Z')^{\sigma'_j}_t = \alpha'_{j_t}(X)$ for $1 \le t \le \ell-1$,
then by the definition of a diagonal basis
(Definition \ref{d:diagbas}, (2)) we have
in the notation of Definition \ref{d:diagbas}
\begin{equation} \label{e:delta} {\rm Res}^{\sigma'_j}
\left (\frac{1}{\alpha'_{k_1}(X) \dots \alpha'_{k_{\ell-1} }(X) } \right )
= \delta_k^j,
\end{equation}
where $\delta_k^j = 1$ if all the multi-indices $k_s = j_s$
(for $s = 1, \dots, \ell-1$) and $\delta_k^j = 0 $ otherwise.
Then we form
$$ Z^{\sigma_j}_{{t}} = (Z')^{\sigma'_j}_t, ~~ 1 \le t \le \ell-1 $$
and
$Z^{\sigma_j}_{\ell} = \hat{e_1}(X). $
(Note that the $Z^{\sigma_j}_\ell$ are all equal for all $j$.)

Then it is clear that  for each $j$ the set which is the union of
$\hat{e_1}$ and the basis $\sigma'_j$
is a basis for $\lies$.
Also, we have that
$${\rm Res}^{\sigma_j}\bigl(\frac{1}{\prod_{t = 1}^\ell Z^{\sigma_k}_t} \bigr ) $$
$$ =
\left ({\rm Res}^{\sigma'_j}   \right )  {\rm Res}_{Z^{\sigma_j}_{\ell} = 0 }
\frac{1}{ Z^{\sigma_j}_{\ell} \prod_{t = 1}^{\ell-1}
(Z')^{\sigma'_k}_t }$$
$$  = {\rm Res}^{\sigma'_j}
\frac{1}{\prod_{t = 1}^{\ell-1} (Z')^{\sigma'_k}_t}  = \delta^k_j
~\mbox{(by (\ref{e:delta}))}. $$
Hence $\sigma$ defined as in the statement of Lemma
\ref{l6.2old} is a diagonal basis for $\bigtriangleup' \cup \{\hat{e_1} \}.$
\hfill $\Box $

Note that in this section we identify weights of $S$ with
weights of $T$ using an invariant inner product on
$\liet$.

\begin{lemma} \label{l6.3old}  We have
$$\emm/\emm_\sigma = \emm'/\emm'_{\sigma'}, $$
where $\emm'$ is the dual lattice of
the  lattice $\enn'$ associated to the hyperplane
arrangement $\bigtriangleup'$ in Lemma \ref{l6.2old},
and $\emm'_{\sigma'}$  is the sublattice of $\emm'$
corresponding to the ordered basis $\sigma'$
for $\lies/\hat{\liet_1}$, and $\sigma = \sigma' \cup
\{\hat{e_1} \}$ is an ordered basis for $\lies$.
Recall that $\emm_\sigma$ was introduced in the statement of
Theorem {\ref{t:szgen}}.
\end{lemma}
\proof By definition $ \emm$ is the
lattice generated by $ \emm'$ and $ \hat{e_1} $
and $\emm'_{\sigma'} $
is the lattice generated by all $\alpha \in \sigma'$
so $\emm_\sigma $ is the
lattice generated by  the $\alpha \in \sigma'$ together with
$ \hat{e_1}. $
Thus the lemma follows immediately.
\hfill $\Box$

Lemmas \ref{l6.2old} and \ref{l6.3old} give the following.

\begin{lemma} \label{l6.4old} Let $\hat{e_1} $ be
a weight of $T$ which generates
the   subgroup $\hat{T_1} \cong S^1$ in $T$.
Here the subgroup $S$ is as in Lemma \ref{l6.2old}.
Then there is a diagonal basis of $T$ consisting of
$\{\hat{e_1} \cup v_{1},  \dots,  \hat{e_1} \cup v_{N} \} $
where $\{v_{1},  \dots,  v_{N} \} $
 is a  diagonal basis for
the hyperplane arrangement in
${\rm Lie} (S/\hat{T_1}) $ corresponding to the weights of the action of
$S/\hat{T_1}$  on the tangent bundle to all components of the
 fixed point set of the action of $S/\hat{T_1}$ on $M_1$.
(Note that here each $v_k$ is a basis for $\lies/\hat{\liet_1}$.)
\end{lemma}

In the following theorem we consider a sum over all
connected subgroups
$S$  of the maximal torus $T$ of $G$  such that $\hat{T}_1 \subset S$, with
Lie algebra ${\bf s}$.
We then consider a sum over
 all connected components $F$ of $M^S$
 for which $T/S$ acts locally freely on $F$
(we refer to these as
$F \in {\mathcal F}(S)_{\rm l.f.}$).

\begin{remark} In Example 7.3 we exhibit an example where
a nontrivial group $T/S$ acts locally
freely on $F \subset M^S$.
\end{remark}

\begin{theorem} \label{t:residgen} \textbf{(Residue formula in the general case)}
Let $\beta$ be an equivariant cohomology class of the
form $\eta e^{\bar{\omega}}$ where $\eta \in H^*_G(M). $
Assume $c \in T$ is a regular value for $\Phi$.

Let $OB$ be the diagonal basis in $\mathcal{OB}(\Delta)$
corresponding to a total ordering of
the elements in
a hyperplane arrangement $\Delta$ $ = \Delta'(S/\hat{T_1}) \cup \{\hat{e}_1\}$
 for ${\bf s}$
(see Lemma \ref{l6.2old}). Here $ \Delta'(S/\hat{T_1})$  is a hyperplane arrangement for
${\bf s}/{\bf t}_1 $
 which comprises
all  weights $\beta_{F,j} $ for the action of $S$ on
the normal bundle to some  component $F$ of the fixed point
set of a connected subgroup $S$ of $T$.
 We sum over all such connected
subgroups $S$.
We have
\begin{equation} \label{e:th6.6}
\int_{\mu^{-1}(0)/G}  \kappa (\eta e^{\bar{\omega}}) =
  \frac{n_1}{n_0'|W|}\int_{\mu^{-1}(\liet)\quott T}
 \kappa (\eta{\mathcal D}  e^{\bar{\omega}})
= \frac{n_1 }{n_0'|W|}
\sum_{S \triangleleft T} \sum_{\sigma \in OB(S) }
\sum_{F \in {\mathcal F}(S)_ {\rm l.f.} }
\end{equation}
$$ ~ \frac{1}{|\emm/\emm_\sigma|} \int_{F\quott (T/S)} \kappa_{T/S}
 \Biggl \{  {\rm Res}^{\sigma} \Bigl (
    {\mathcal D}^2(X)
\sum_{m \in \emm/\emm_\sigma}
 \frac{\eta(X) e^{(\mu(F)  -m)(X) } e^{\omega_F} }
{ e_F(X) \prod_{\alpha \in \sigma}
( 1 - e^{\alpha(X) } ) } \Bigr ) ~ \Biggr \}, $$
where the notation ${\rm Res}^\sigma$ was introduced in
Definition \ref{e:defres}.
Here $X \in \liet \times {\mathbb C}$.
The notation $e_F$ refers to the $S$-equivariant
Euler class for the action of $S$ on the normal bundle to $F$
in $M$.
Also,
$n_1$  and $n_0'$
were introduced in Proposition \ref{pr:redab}.
\end{theorem}

\proof
We shall perform an  induction reducing the rank of the group
$T$ that acts on $M$, replacing $M$ by a component
$M_1$ of $M^{\hat{T_1}}$,  replacing $T$ by $T/\hat{T_1}$
and replacing $\eta$ by ${\rm Res}_{1} \eta$ in the notation of
(\ref{e:defres}). Finally we reach a
space $M_i$ and a group $T_i$ for which $T_i$
acts locally freely on $M_i$.

Let $\eta  \in H^*_T(\widetilde{M})$;
then
by
 (\ref{6.1})  combined with periodicity
(in other words with Proposition  \ref{p:perone})  we have
\begin{equation} \label{e:stepone}
\int_{\mu^{-1}(\liet)\quott T }
\kappa_T ({\mathcal D} \eta e^{\bar{\omega}} )
\end{equation}
\[
= \sum_{\widetilde{M}_1}\int_{\mu^{-1}(\liet) \cap{\widetilde{M}_1} \cap
\mu_{T/\hat{T}_1}^{-1}(0)/(T/\hat{T}_1) }
\kappa_{(T/\hat{T}_1)} {\rm Res}_{Y_1 = 0 }
\left ( \frac{ {\mathcal D} \eta e^{\bar{\omega}}  }{(1-e^{Y_1})
e(\tilde{\nu})   } \right ).
\]
Here
 $\hat{e_1} \in \liet$ generates $\hat{T_1}$ and
$\widetilde{M_1}$ is a component of $\widetilde{M}^{\hat{T}_1}$
for which $0 \le \langle \mu_{\hat{T}_1}(\widetilde{M}_1),
\hat{e_1} \rangle\linebreak \le \langle \hat{e_1}, \hat{e_1}
\rangle  $ (the set of such $\widetilde{M_1}$ is in bijective
correspondence with the components $M_1$ of $M^{\hat T_1}$). Also
$\tilde{\nu}$ is the normal bundle to $\widetilde{M_1}$ in
$\mu^{-1}(\liet)$ and $e(\tilde{\nu})$ is its
 equivariant Euler class in $\mu^{-1}(\liet)$.
These objects are in bijective correspondence with the
corresponding objects $\nu$ and $e(\nu)$ in $M$, because
when we restrict to $\mu^{-1}(\liet)$, the exponential map
becomes a local  diffeomorphism so the
projection map from
$\widetilde M \cap \mu^{-1}(\liet)$ to
$M \cap \Phi^{-1}(T)$  is a covering map.

Recall from Remark \ref{r:slice} that the restriction of
$\mu_T$ to $\mu^{-1}(\liet)$ is proper.
Hence Proposition \ref{p:perone} can be applied.
In the first paragraph of the proof of Lemma
\ref{l:4.1} (referring to Lemma \ref{l3.3a})
we establish the existence of regular values so that
we can apply Proposition \ref{p:perone}.

We assume by induction on $r$  that
\begin{equation} \label{e:steptwo}
\int_{ \mu^{-1}(\liet))^{\hat{T}_1} \quott (T/\hat{T}_1) }
\kappa_{r-1} ({\mathcal D}\eta e^{\bar{\omega} })
 =
\sum_{S' \triangleleft  T/{\hat{T}_1 }}
\sum_{\sigma' \in OB(\Delta'(S')) }
 \sum_{F' \in {\mathcal F}(S')_{\rm l.f.}}
\end{equation}
$$
\frac{1}{|\emm'/\emm'_{\sigma'}|}{\rm Res}^{\sigma'}  \int_{F'\quott (T/S)}
\kappa_{T/S}
 \left \{ \eta(X') {\mathcal D} (X')
\sum_{m' \in \emm'/\emm'_{\sigma'}} \frac{e^{ (\mu(F')-m')(X') }e^{\omega_{F'} } }
{e_{F'}  \prod_{\alpha' \in \sigma'}
( 1 - e^{\alpha'(X') } )} \right \}. $$
where $OB(\Delta'(S'))$
is a diagonal  basis  associated to a hyperplane arrangement
 $\Delta'(S') $ in $\liet/\hat{\liet}_1$
consisting of
 all weights for the action of $S'$ on
the normal bundle to $F'$ for all $F' \in {\mathcal F} (S')_{\rm l.f.}$
(see Lemma \ref{l6.2old})
and $e_{F'}$ is the $(S' := S/\hat{T}_1)$-equivariant Euler class of the normal
bundle in $\mu^{-1}(\liet)^{\hat{T_1}}
$ to a component $F'$ of the fixed point set of the
action of $S'$ on $\mu^{-1}(\liet)$, while $X' \in \liet/\hat{\liet}_1
\otimes \CC$.
Notice that $F' \in {\mathcal F} (S')$ inherits a Hamiltonian
action of $T/S'.$ If $T/S'$ acts locally freely on $F'$,
then we cannot proceed further with the induction.
Otherwise, continuing the induction we finally replace $F'$ by the fixed point set
$F'' \subset F'$ of a subgroup
$S''$ for which $T/S''$ acts locally freely on
$F''$.

\begin{remark} \label{r:bij}
Note that if we define $S$ to be the pullback of the subgroup
$S'$ under the
projection map $T \to T/\hat{T_1}$, it is clear that
$T/S \cong T'/S'$ where $T' = T/{\hat{T_1}}.$
\end{remark}

Combining (\ref{e:stepone}) with (\ref{e:steptwo})
we obtain
\begin{equation} \label{e:stepthree}
\int_{\mu^{-1}(\liet)\quott T} \kappa ({\mathcal D} \eta
e^{\bar{\omega}} ) = \sum_{\widetilde{M_1}} \sum_{S' \triangleleft
T/\hat{T_1}} \sum_{\sigma' \in OB(\bigtriangleup'(S')) } \sum_{F'
\in {\mathcal F} (S')_{\rm l.f.} }
\frac{1}{|\emm'/\emm'_{\sigma'}|} {\rm Res}^{\sigma'} {\rm
Res}_{Y_1 = 0 } \end{equation}
$$\int_{F' \quott (T/S)} \kappa_{T/S} \Bigl \{
\eta(X') {\mathcal D}_H (X') \sum_{m' \in \emm'/\emm'_{\sigma'}}
\frac{e^{(\mu(F') - m) (X')} e^{\omega_{F'} } } {e_{F'}
e_{\widetilde{M_1}}(1  - e^{Y_1}) \prod_{\alpha' \in \sigma'} (1 -
e^{\alpha'(X')}) } \Bigr \}.
$$
By Lemma  \ref{l6.2old} there is a bijective correspondence
between the elements $\sigma'$ in the diagonal
basis  $OB(\bigtriangleup'(S'))$ and
the elements $\sigma$ in a diagonal basis
  $OB(\bigtriangleup(S))$.

Note that
$e_{F'} e_{\widetilde{M_1}} $ is the equivariant
Euler class  $ e_F^\liet$ of the normal bundle to $F$ in
$\mu^{-1}(\liet)$.
The equivariant Euler class of the normal bundle to $F$ in $M$ is
$e_F = e_F^\liet {\mathcal D}$.
Collecting these results  we obtain  (\ref{e:th6.6}).

By Lemma   \ref{l6.2old} the denominator of (\ref{e:stepthree}) is
$e_F \prod_{\alpha \in \sigma}
(1  - e^{\alpha(X)} )$ for $\sigma$ the element
of
  $OB(\bigtriangleup(S))$ formed from $\sigma'$ and $\hat{e_1}$.
This lemma also exhibits a bijective correspondence
between
  $OB(\bigtriangleup'(S'))$ and
  $OB(\bigtriangleup(S))$.
By Lemma \ref{l6.3old}
$|\emm'/\emm'_{\sigma'}| =
|\emm/\emm_{\sigma}|$.
In Remark \ref{r:bij}
we have exhibited a bijective correspondence
between connected subgroups
$S$ containing $\hat{T_1}$ and
connected subgroups
$S'$ of $T/\hat{T_1}$.
The space $F'$ is a component of
$\Bigl ( \mu^{-1}(\liet)^{\hat{T_1}}\Bigr )^{S'}$ on
which $T'/S'$ acts locally freely, so it is  a component
$F' =F$ of $\mu^{-1}(\liet)^S$ on which $T/S$ acts
locally freely.
Finally the residue
${\rm Res}^{\sigma'} {\rm Res}_{Y_1 = 0 }$
is equal to $
{\rm Res}^{\sigma}. $

This
gives the result, using Lemmas \ref{l6.2old},  \ref{l6.3old} and
\ref{l6.4old} to identify the diagonal bases.

Note that
the finite sum
$$ \sum_{m \in
\emm/\emm_{\sigma}}
\frac{e^{-m(X)} }{
\prod_{\alpha \in \sigma}
( 1 - e^{\alpha(X) } )} $$
parametrizes a class of  different components of the fixed point
set  of the $T$ action
in $\widetilde{M} \cap \mu^{-1} (\liet)$ which correspond
    to a given $F \subset M\cap \Phi^{-1}(T)$
(by specifying  different values of
the moment map $\mu$ which  correspond to the same
value of $\Phi(F)$).

\hfill $\Box$

 \section{Alekseev-Meinrenken-Woodward}

  As we explained above,   in
 \cite{AMW1} Alekseev,  Meinrenken and  Woodward
 gave a  formula
for intersection pairings of   reduced spaces of a
 quasi-Hamiltonian $G$-space. In this section,
we will relate our formula with Alekseev-Meinrenken-Woodward's by using
 Szenes' theorem (Theorems \ref{t:szenes} and  \ref{t:szgen}).

 Let $\Lambda^w$ denote
 the weight lattice in $\liets$.
 For $\lambda \in \Lambda^w$
 let $\chi_{\lambda}$ denote the character of the corresponding irreducible representation $V_{\lambda}$.

Let $c \in Z(G)$.
If $c$ is a regular value of $\Phi$ let
\[
     \kappa : { H}^*_G(M) \rightarrow H^*(M_{\rm red})
\]
be the composition of
 the pull-back to the level set
$H_G^*(M) \rightarrow H_G^*(\Phi^{-1}(c))$ and
the isomorphism $H_G^*(\Phi^{-1}(c)) \cong H^*(\Phi^{-1}(c)/G)$.
(Here by analogy with
(\ref{e:reddef}) we define $M_{\rm red} = \Phi^{-1}(c)/G$,
cf. Theorem \ref{t:qhr}.)
 For any $\lambda \in \Lambda^w$, all fixed point manifolds $F \in {\mathcal F}(\lambda + \rho)$ (defined in the statement of Theorem \ref{t7.5})
 are contained in the preimage of the maximal torus $\Phi^{-1}(T)$. This follows by equivariance
 of $\Phi$.
 In fact, $F$
 is a quasi-Hamiltonian $T$ space in the fixed point set of $S_1$ (where
 $S_1$  is the connected subgroup of $G$ generated by $\lambda + \rho$),
 with $\omega_F = \iota_F^* \omega$ as its symplectic form and $\Phi_F
 = \Phi |_F$ as its moment map. Since $\Phi_F$ takes its value in $T$, we can compose with the map $T
 \rightarrow U(1), t \mapsto t^{\lambda + \rho}$ and this composition is denoted $(\Phi_F)^{\lambda + \rho}$.
 Note that $(\Phi_F)^{\lambda + \rho}$ is constant along $F$.

\begin{theorem}[Alekseev-Meinrenken-Woodward \cite{AMW1}] \label{t7.5}

   Let $G$ be the direct product of a connected, simply connected Lie group
and a torus, and
   ($M$,$\omega$,$\Phi$) a compact quasi-Hamiltonian $G$-space.
Suppose
 $\beta$ is a class of the form
$\beta =\eta e^{ \bar{\omega}} $
where $\eta \in H^*_G(M)$  and
$\bar{\omega}$ was defined in (\ref{e:barwdef}).
Then we have (\cite{AMW1},
first equation p. 346, after the  proof of Theorem 5.2)
 \begin{equation} \label{eq:amw}
    \int_{M_{\rm red}} \kappa(\eta) \exp(\omega_{red})
 = \end{equation}$$ \frac {k} { (vol G)^2} \sum_{\lambda \in \Lambda^w_+}
  \left( d_\lambda^2
 \sum_{F \in {\mathcal F}  (\lambda + \rho)}
 (\Phi_F)^{\lambda + \rho} \int_F \frac {\iota_F^* \eta (2 \pi i (\lambda + \rho))}
 {e_F( 2 \pi i (\lambda + \rho))} e^{\omega_F} \right)$$

 where $k$ is the order of the principal stabilizer
and $d_\lambda = {\rm dim} ~V_\lambda$.
Here, ${\mathcal F} (\lambda + \rho)$ is the set of  components
of the zero locus of the vector field associated to
$\lambda + \rho$.
We use the notation $\Lambda^w_+$ to denote the set of weights in
the fundamental Weyl chamber, i.e. the dominant weights.
 \end{theorem}

Note that in \cite{AMW1}, the sum over $\lambda\in \Lambda^w_+$
is actually over $\Lambda^w_+\subset \liet_+$, the image
in $\liet$ of the intersection
of the fundamental Weyl chamber  in $\liet^*$ with the weight lattice
 under the bijective map  $\liet^* \to  \liet$ given
by the bi-invariant inner product. This map identifies
the integer lattice in $\liet$  with a sublattice
of the image of the weight lattice.

To study the formula in the above theorem we need some tools from
Lie  group  theory.

\begin{theorem} [Weyl dimension formula, \cite{Lie}]
   If the irreducible representation $V_{\lambda}$  of $G$ has highest weight $\lambda$, then its dimension is given by
  \[
       dim V_{\lambda} = \prod_{\alpha \in R_+(G)} \frac {
\langle\alpha, \lambda + \rho \rangle} {\langle \alpha, \rho\rangle},
  \]
  where $\langle \cdot,\cdot\rangle$
is the  $W$-invariant inner product on ${\bf t}$ corresponding
to the bi-invariant inner product on $\bf g$, and
$R_+(G)$ are the positive roots of $G$.
\end{theorem}

\noindent{\bf Comparing  Theorem \ref{t:residgen} with Theorem \ref{t7.5}:}
We will compare the formula (\ref{eq:amw}) from \cite{AMW1}
with our formula (\ref{e:th6.6})
in Theorem \ref{t:residgen},
when $G$ is a compact Lie group
which is the product of  a connected simply-connected
Lie group with a torus.

\begin{remark}In fact Theorem \ref{t7.5} is
valid for arbitrary compact connected $G$ if the formula is interpreted
in an appropriate way, using a finite cover $G'$ of $G$ which is
a product of a connected simply-connected group and a torus and
an associated quasi-Hamiltonian $G'$-space $M'$ which is a finite cover
of $M$. For simplicity the formula in \cite{AMW1} Theorem 5.2
 is stated under the hypothesis that
$G$ is the product of a connected simply connected Lie group and
a torus. See Section 5.2 of \cite{AMW1}.

Our proof invokes the hypothesis that $G$ is the product of a connected
simply-connected group and a torus to identify
$\rho$ (half the sum of the positive roots of $G$) with an element in
the  weight lattice of $G$
 (see Remark \ref{r:meinhyp}).

\end{remark}

We will use Szenes' formula (Theorem
\ref{t:szgen}) to relate (\ref{eq:amw}) to our formula.
We recast  (\ref{eq:amw}) (replacing $2 \pi i (\lambda + \rho)$
by $\lambda + \rho$) as (\ref{e:sevenone}) below.
We observe that
$T$ acts on  each component $F_1$ of the fixed point set of $\lambda + \rho$,
so we can apply the Atiyah-Bott-Berline-Vergne localization theorem to replace
$$ \int_{F_1} \frac{ \eta(\lambda + \rho)}{e_{F_1} (\lambda + \rho) } $$
by
$$ \sum_{F \subset (F_1)^S} \int_{F} \frac{ \eta(X)}{e_{F} (X) }|_{X = \lambda + \rho} $$
for any subgroup $S$ of $T$ containing the  subgroup
$S_1$ generated by $\lambda + \rho$ (using the fact that the equivariant
Euler class is multiplicative so $e_F$ is the product of
$e_{F_1}$ with the equivariant Euler class of the normal bundle to
$F_1$ in $F$). We replace $F_1$ by $F$ and
$S_1$ by $S$ and continue inductively
until we reach a  subgroup $S$ for which
$T/S$ acts locally freely on $F$.

\begin{remark} \label{r:meinhyp}
 Recall  that
since  $G$ is assumed to be the product of a connected simply connected group
and a torus, $\rho \in \Lambda^w$.
So we replace the sum over
$\lambda + \rho$ by a sum over $\xi \in \Lambda^w_0 \subset \liet$,
where $\Lambda^w_0$ is the set of strictly
dominant weights, in other words those weights in $\liet_+$
not on the boundary
of $\liet_+$ (not on any hyperplanes specified by one of the fundamental
weights). We have  used Weyl invariance of the sum to replace the sum
over strictly dominant weights by a sum over the regular weights
(those weights in $\liet$ not on any hyperplanes specified by
fundamental weights).

\end{remark}

Hence we obtain

\begin{equation} \label{e:sevenone}
    \int_{M_{\rm red}} \kappa (\eta e^{\bar{\omega}})
  =   \frac{k} { |W| (VolG)^2}
\sum_{S \triangleleft T} \sum_{\lambda + \rho \in \Lambda^w_0(S) }
     (dimV_\lambda)^2 \end{equation}$$\sum_{F \subset (M^S)_{\rm l.f.}}
(\Phi_F)^{\lambda + \rho}
     \int_F \frac {\eta(\lambda + \rho)} {e_F(\lambda + \rho)} e^{\omega_F},
$$
where
$M^S$ is the fixed point set of the action of $S$ on $M$
and $e_F$
is the equivariant Euler class for the normal bundle to $F$ in $M$.
Here $(M^S)_{\rm l.f.}$ denotes the components $F$ of the fixed point
set of $S$ acting on $M$ for which $T/S$ acts locally freely  on $F$;
this is indexed by the set ${\mathcal F}(S)_{\rm l. f.}$.
We subdivide  $\Lambda^w_0$ into
collections   $\Lambda^w_0(S)$ of regular
weights $\xi$ in the Lie algebra of   a  connected subgroup $S$
(those regular weights of $G$ which lie in $\lies$).

Now the formula (\ref{e:sevenone})
becomes

\begin{equation} \label{e:7.4}
    \int_{M_{\rm red}} \kappa (\eta e^{\bar{\omega}})  =     \frac{k}
{|W| (VolG)^2}
\sum_{S  \triangleleft T}
     \sum_{ \xi \in  \ \Lambda^{w}_0(S) } \sum_{F \in  {\mathcal F}(S)_{\rm l.f.}  }
     \end{equation}$$\prod_{\alpha \in R_+(G)} \frac {\langle  \alpha, \xi \rangle ^2} { \langle \alpha, \rho \rangle ^2}
     e^{ \langle  \mu(F), \xi \rangle }
     \int_F \frac {\eta(\xi)} {e_F(\xi)} e^{\omega_F}.
$$
(Here, for each $F$ we have chosen $\mu(F) \in \liet$
for which $\Phi(F) = \exp \mu(F).$
The quantity $     e^{ \langle  \mu(F), \xi \rangle } $
is well defined because $\xi \in \Lambda^w.)$
Now we apply Szenes' theorem (separately for each $S$) to
\begin{equation} \label{e:6.3'}
{\mathcal S}_S :=       \sum_{ \xi \in \Lambda^w_0(S) }
\sum_{F \in  {\mathcal F}(S)_{\rm l.f.}}
  e^{ \langle  \mu(F), \xi \rangle } f_F(\xi)
\end{equation}
where
\[
     f_F(X) = \prod_{\alpha \in R_+(G)} \frac {\langle  \alpha, X \rangle ^2}
       { \langle \alpha, \rho \rangle ^2}  \int_F  \frac{\eta(X) e^{\omega_F}} {e_F(X)}
\]
for $X \in {\rm Lie}(T)$.

According to Szenes' theorem (Theorem \ref{t:szgen}),
 \[  {\mathcal S}_S =
    \sum_{\sigma \in OB(\Delta^R(S) )} \sum_{F \in {\mathcal F}(S)_{\rm l.f.}}
 \]\[{\rm Res}^\sigma \left( e^{ \langle \mu(F), X_S \rangle }
                 f_F(X_S)
            \frac{1}{|\emm/\emm_\sigma|} \sum_{m \in \emm/\emm_\sigma}
              \frac {e^{-m(X_S)}}
{\prod_{\alpha \in \sigma} (1-e^{\alpha(X_S)})}
 \right).\]

Here $X_S\in {\rm Lie}(S)
 \otimes \CC$ and  $\emm\subset \liet^* $ is the
dual lattice
to the
lattice $\Lambda^w(S)$.
We use the hyperplane arrangement $\Delta^R(S)$  given by the
restriction of the  simple roots of $G$ to the
Lie algebra  ${\lies}$ of $S$,
and choose a diagonal basis $OB(\Delta^R(S))$ for
this hyperplane arrangement.

This gives for the final form of  Theorem \ref{t7.5}:
$$
    \int_{M_{\rm red}} \kappa(\beta) \exp(\omega_{red})   =
\frac{k} {|W| (VolG)^2}  \sum_{S \triangleleft T} {\mathcal S}_S
$$
\begin{equation} \label{e:7.5}  = \frac{k} {|W| (VolG)^2}
    (\prod_{\alpha \in R_+(G)} \frac {1} { \langle \alpha, \rho \rangle ^2})
\sum_{S\triangleleft T}
 \sum_{F \in {\mathcal F}(S)_{\rm l.f.} }
   \sum_{\sigma \in OB(\Delta^R(S))}
\end{equation}
$$ \frac{1}{|\emm/\emm_\sigma|} \Bigl ( {\rm Res}^\sigma  \
 {\mathcal D}^2(X_S)
   \sum_{m \in \emm/\emm_\sigma}
\frac{     e^{(\mu(F)-m) (X_S)}}{\prod_{\alpha \in \sigma} (1-e^{\alpha(X_S)})}
\int_F     \eta(X) \frac { e^{\omega_F} }
                          {e_F(X_S) }
   \Bigr  ) ,  $$
where the notation ${\rm Res}^\sigma$ was introduced in
Definition \ref{d:diagbas}. Here $X_S$ is a variable in $\lies \otimes \CC$.

Comparing (\ref{e:7.5}) to our formula (\ref{e:th6.6})
in Theorem \ref{t:residgen}
the differences are:

\noindent{(1)} The presence in Theorem \ref{t:residgen}
of $\kappa_{T/S}$ and evaluation on
$F\quott (T/S)$ rather than $F$ as in Theorem \ref{t7.5}.
The identification between the   two formulas follows because the
following result  is true:
\begin{prop} \label{p:amwspcase}
If a locally free
action of $T/S$    on a  quasi-Hamiltonian $T/S$-space $(F, \Phi)$
induces an action of $T$ on $F$,
 then  for all $\eta \in H^*_T(F)$
and corresponding $\beta = \eta e^{\bar{\omega}} $
\begin{equation} \label{e:conj}
\int_F \beta \circ i_S^* = \int_{F\quott (T/S)} \kappa_{T/S} \beta
\end{equation}
where the inclusion $i_S: S \to T$ induces
$i_S^* : H^*_T(F) \to H^*_S(F) $
(setting the variables $X $ corresponding to $T/S$ to
zero)
and both sides of (\ref{e:conj})
are in $H^*_S$, and
$\kappa_{T/S} : H^*_T(F) \to
H^*_S(\Phi^{-1}(e)/(T/S))$
is the Kirwan map (in other words the restriction map
to $\Phi^{-1}(e)$ composed with the isomorphism
$H^*_{T/S} (\Phi^{-1}(e)) \cong H^*(\Phi^{-1}(e)/(T/S))$).
\end{prop}

To prove this result  it suffices to consider the case $S = \{e\}$
for which it reduces to verifying that
\begin{equation} \label{e:jk}
\int_{\Phi^{-1}(e) \cap F/T}  e^\omega  \eta (X_j = c^{(j)}) =
\int_F  e^\omega \eta|_{ (X = 0)}.
\end{equation}
(Here  $c^{(j}$ is the first Chern class of the line bundle
associated to the principal $T$-bundle
\begin{equation} \label{e:bun}
T \to \Phi^{-1}(e) \mapsto \Phi^{-1}(e)/T \end{equation}
by the $j$-th weight $T \to U(1)$.)

Proposition \ref{p:amwspcase}
 is clear when the $c^{(j)}$ are $0$ (i.e. when the bundle
(\ref{e:bun}) is trivial). In general it is
a special case of Theorem 5.2 of \cite{AMW1}
(the case where the group $G$ is a torus and acts locally
freely on $M$). We expect it may  be possible to give a
direct proof of this special case without relying on the
full strength of Theorem 5.2 of \cite{AMW1}.
  For a related result, see \cite{jw}.

\noindent{(2)}
Differences between overall multiplicative constants
\begin{equation} \label{constant3}
     \frac{k} {|W| (VolG)^2}
    (\prod_{\alpha \in R_+} \frac {1} { \langle \alpha, \rho \rangle ^2}).
\end{equation}
and
\[
     \frac{n_1}{|W|n_0'}
\]
Since  by \cite{BGV} Cor. 7.27
\begin{equation} \label{e:bgvvol}
\frac {vol G}{vol T} = \prod_{\alpha \in R_+(G)} \frac{1}{2\pi \langle \alpha, \rho \rangle}
\end{equation}
where the volumes are taken with respect to the metrics
on $G$ and $T$ given by the bi-invariant inner product $\langle \cdot , \cdot
\rangle$,
the equation (\ref{constant3}) is equal to
\[
   \frac { (2 \pi)^{2n_+}  k} { |W| (vol T )^2}.
\]
So the difference
 is $\frac{n_1}{|W|n_0'
}$ versus $\frac{(2 \pi)^{2n_+}  k} {|W| (vol T)^2}$.
This difference arises from the
normalization of the Riemannian metrics used in  (\ref{e:bgvvol}).

\noindent{(3)} The hyperplane arrangement
$\Delta^R(S)$ used in (\ref{e:7.5}) is different from the
hyperplane arrangement $\Delta(S) $ used in Theorem \ref{t:residgen}.
If the weights at all components $F \in {\mathcal F}(S)$
 are roots,
the same hyperplane arrangement can be used. This is true for example
in the case  studied in \cite{JK3}.

\section{Applications}

\noindent{\bf Example 7.1:}
The main example for which the residue formula
has been used previously to compute intersection numbers
in reduced spaces of quasi-Hamiltonian spaces appears in  \cite{JK3}, for which
$M =
G^{2h}$ (where $G=SU(n)$ and $G$ acts by conjugation on $M$). In this example
the fixed point sets of all vector fields $v_\lambda$ associated to
weights $\lambda$ are equal to $M^T = T^{2h}$.

\noindent{\bf Example 7.2:} If $T$ is a torus, a  standard example of a
quasi-Hamiltonian $T$-space is $M = T \times T$, with the action of $T$ given
by multiplication of the first copy of $T$. In this case $T$ acts freely
on $M$.
The moment map is the projection to the second factor $T$.
The reduced space for this action is a point.

\noindent {\bf Example 7.3:} An example due to C. Woodward \cite{Wood}
occurs when $G = SU(3)$, with rank two maximal torus $T $ for which $M
$ is constructed starting with $ G \times_T \CC P^2 $ by replacing
three orbits $G/T$ (corresponding to the three fixed points of $T$ in
$\CC P^2$) by $G/U(1)_j$ (for $j = 0, 1,2$).
 Here $U(1)_0\subset T$ is generated by
$\rho$
(where $\rho$ is  half the sum of the positive roots, identified as
an element of
${\bf t}$ using the bi-invariant inner product) and
$U(1)_j$ is generated by $c_j \rho$ (for $j = 1, 2$) where $c_j $
are nontrivial elements of the Weyl group (the cyclic permutations).
The
fixed point set $F_j$ of $U(1)_j$ acting on $M$ is symplectic, and is
isomorphic to $S^1 \times S^1$ with a multiple of the standard
symplectic structure
(this is the case treated in
Example 7.2).  The normal bundle of $F_j$ is trivial, and is
isomorphic (as a representation of $U(1)_j$) to $\CC \oplus
\liet^\perp $ (the direct sum of the trivial representation and the
restriction from $T$ to $U(1)_j$ of the adjoint representation of $T$
on the orthocomplement $\liet^\perp$ of $\liet$ in $\lieg$,
where $U(1)_j$ acts on $\liet^\perp$ via   $c_j \rho$.)  Thus the
Euler class is
$$ e_{F_j} (Y_j) = \widetilde{\gamma_1}(Y_j) \widetilde{\gamma_2}(Y_j)
(\widetilde{\gamma_1}(Y_j) + \widetilde{\gamma_2}(Y_j) )$$
where $\gamma_1$ and $\gamma_2$ are the simple  roots of
$SU(3)$ and the third positive root is $\gamma_1 + \gamma_2$, and
we have denoted by $\widetilde{\gamma_j}$ the restrictions of $\gamma_j$
to the Lie algebra of $U(1)_j$. Here we have introduced
a nonzero element $Y_j$ of the Lie algebra of $U(1)_j$.

In this example $M^T = \emptyset$.
The orbit type decomposition consists of two connected strata:
the principal stratum where $G$ acts freely, and another
stratum $M_1$ where the stabilizer is conjugate to $U(1)$.

By Theorem \ref{t:residgen} the relevant subgroups $S$ are the
$U(1)_j$ conjugate to $U(1)$; the group $T/S$ is also isomorphic to
$U(1)$, and $T/S$ acts freely on the fixed point set $F_j$ of
$U(1)_j$. So there is only one element in the hyperplane arrangement
(which we denote by  $\alpha_j$)
and $\emm_\sigma = \emm$.  So the expression in Theorem
\ref{t:residgen} reduces to
\begin{equation}
\sum_{j=0}^2 \int_{F_j \quott (T/U(1)_j)} \kappa_{T/U(1)_j}
\Biggl \{ {\rm Res}_{Y_j = 0 }
~ \Bigl ( \frac{ {\mathcal D}^2 (X) \eta(X) e^{\mu(F_j)(X)} e^{\omega_{F_j}} }
{ e_{F_j} (1 - e^{\alpha_j (X)} ) } \Bigr )~ \Biggr \}
\end{equation}
where $Y_j \in {\rm Lie} (U(1)_j)$ is as above.
In this example, $F_j \cong U(1) \times U(1) $ so $F_j \quott (T/U(1)_j)
$ is a point, as in Example 8.2.
Thus  the final
map $\kappa_{T/U(1)_j} $ in our computation is
$$\kappa_{U(1)}: H^*_{U(1) } (U(1) \times U(1))
\cong H^*(U(1)) \to H^*({\rm point}) = H^*(\Phi^{-1}(e)/U(1)),$$
which is the pullback map.

\smallskip

\noindent{\bf Acknowledgments:} This article presents
the main results of the Ph.D. thesis \cite{song} of the
second author, written under the supervision of the
first author. The authors would like to acknowledge the
helpful comments of Rebecca Goldin, Yael Karshon
and Eckhard Meinrenken (who formed the committee of examiners for
the second author's Ph.D. thesis).
We wish to thank Anton Alekseev for pointing out that in an
earlier version of Section 6 we had overlooked the distinction
between a sum over the entire weight lattice and the dominant weights;
the latter was used in  Theorem 5.2 of \cite{AMW1}.
We also wish to thank the referee for useful
comments on an earlier version, and for alerting us to some
relevant results
in the literature.

\noindent Lisa Jeffrey  and  Joon-Hyeok Song \\ Department of
Mathematics\\ University of Toronto\\ Toronto, ON M5S 2E4, Canada\\
E-mail: jeffrey@math.toronto.edu\\
song@math.toronto.edu


\begin{thebibliography}{99}

  \bibitem{AMM} A. Alekseev, A. Malkin, E. Meinrenken, Lie group valued moment maps
          {\em J. Differential Geometry} {\bf 48} (1998) 445-449.

  \bibitem{AMW1} A. Alekseev,
E. Meinrenken, C. Woodward, Group-valued equivariant localization.
{\em Inventiones Mathematicae} {\bf 140} (2000) 327-350.

  \bibitem{AMW2} A. Alekseev,
E. Meinrenken, C. Woodward, Duistermaat-Heckman measures and moduli
          spaces of flat bundles over surfaces.
          {\em Geometric and Functional Analysis} {\bf 12} (2002) 1--31.

\bibitem{AMW3}  A. Alekseev,
E. Meinrenken, C. Woodward, The Verlinde formulas as fixed point
formulas, {\em J. Sym. Geo.}{\bf 1} (2001) 1-46.


  \bibitem{AB} M.F. Atiyah, R.Bott, The moment map and equivariant cohomology.
{\em Topology} {\bf 23} (1984) 1-28.

\bibitem{BGV} N. Berline, E. Getzler, M. Vergne,
{\em Heat Kernels and Dirac operators} (Grundlehren \#298),
Springer-Verlag, 1992.

  \bibitem{BV} N. Berline, M. Vergne,
Z{\'e}ros  d'un champ de vecteurs et classes
          caract\'eristiques  \'{e}quivariantes.
          {\em Duke Math. J.} {\bf 50} (1983) 539-549.
\bibitem{BV1} M. Brion, M. Vergne, Arrangement of hyperplanes I: Rational
functions and Jeffrey-Kirwan residue. {\em Ann. Scient. Ec. Norm.
Sup.} {\bf 32} (1999) 715-741.

\bibitem{BV2} M. Brion, M. Vergne, Arrangement of hyperplanes II:
Szenes formula and Eisenstein series. {\em Duke Math. J.}
{\bf  103} (2000) 279-302.



  \bibitem{Lie} J. J. Duistermaat, J.A.C. Kolk, {\em Lie Groups.}
Springer-Verlag, 2002.

\bibitem{Goldin}
 R.F. Goldin,  An effective algorithm for
the cohomology ring of symplectic reductions,  {\em Geom. Anal.
Funct. Anal.} {\bf 12}
 (2002) 567--583.


  \bibitem{GK} V. Guillemin, J. Kalkman, The Jeffrey-Kirwan localization theorem and
          residue operations in equivariant cohomology.
{\em J. reine angew. Math.} {\bf 470} (1996) 123-142.





\bibitem{J} L.C. Jeffrey {\em Extended moduli spaces of flat connections on Riemann surfaces}
Math. Annalen {\bf 298} (1994) 667-692.

  \bibitem{JK1} L.C. Jeffrey, F.C. Kirwan, Localization for nonabelian group actions.
{\em Topology} {\bf 34} (1995) 291-327.





  \bibitem{JK3} L.C. Jeffrey, F.C. Kirwan,
Intersection theory on moduli spaces
          of holomorphic bundles of arbitrary rank
on a Riemann surface. {\em Annals of
          Mathematics} {\bf 148} (1998) 109-196.

  \bibitem{K} F. Kirwan, {\em Cohomology of quotients in algebraic and symplectic geometry.}
          {Math. Notes vol. 31, Princeton Univ. Press, Princeton, NJ, 1985}



  \bibitem{martin} S.K. Martin,
 Symplectic quotients by a nonabelian group and by its maximal torus,
preprint math.SG/0001002;
Transversality theory, cobordisms, and invariants of symplectic quotients,
preprint math.SG/0001001; {\em Ann. Math.}, to appear.

\bibitem{Wood} E. Meinrenken and C. Woodward,
 {\em Chris Woodward's example}, unpublished manuscript (2001).


\bibitem{plamen} O. Plamenevskaya, A residue formula for
$SU(2)$-valued moment maps,
preprint math.DG/9906093; {\em Can.  Bull. Math.},  to appear.

\bibitem{song} Joon-Hyeok Song, {\em Intersection numbers in
reduced spaces of q-Hamiltonian $G$-spaces}, Ph.D. thesis,
University of Toronto, 2004.

\bibitem{Sz} A. Szenes, Iterated residues and multiple Bernoulli
polynomials, {\em International Mathematics
Research Notices} {\bf 18} (1998) 937-956.



\bibitem{jw} J. Weitsman, A Duistermaat-Heckman formula for symplectic
circle actions, {\em International Mathematics
Research Notices} {\bf 12} (1993) 309-312.






 \end{thebibliography}
 \end{document}